\documentclass[aip,amsmath,amssymb,reprint]{revtex4-2}

\usepackage{graphicx}
\usepackage{dcolumn}
\usepackage{bm}

\usepackage{color}
\usepackage[utf8]{inputenc}
\usepackage[T1]{fontenc}
\usepackage{mathptmx}
\usepackage{datetime}

\newcommand{\ie}{\emph{i.e.}}
\newcommand{\eg}{\emph{e.g.}}
\newcommand{\norm}[1]{\left\Vert#1\right\Vert}

\newcommand{\SM}{S}
\newcommand{\me}{m}

\begin{document}


\title{How Can We Describe Density Evolution Under Delayed Dynamics?}
\author{Michael C. Mackey}
\email{michael.mackey@mcgill.ca}
\homepage{https://www.mcgill.ca/mathematical-physiology-lab/}
\affiliation{Departments of Physiology, Physics \& Mathematics, McGill University,  Montreal, Quebec, CANADA, H4X 2C1}
\author{Marta Tyran-Kami\'{n}ska}
\email{ mtyran@us.edu.pl}
\affiliation{Institute of Mathematics,
University of Silesia,
Bankowa 14,
40-007 Katowice,
Poland}

\begin{abstract}
Although the theory of density evolution in maps and ordinary differential equations is well developed, the situation is far from satisfactory in continuous time systems with delay. This paper reviews some of the work that has been done numerically and the interesting dynamics that have emerged, and the largely unsuccessful attempts that have been made to analytically treat the evolution of densities in differential delay equations. We also present a new approach to the problem and illustrate it with a simple example.
\end{abstract}
\maketitle

\begin{quotation}
In this paper we highlight an open problem in mathematics that has implications for any system whose dynamics are dependent on behavior in the past.  Namely how can we describe the evolution of {\it densities} in such systems. We review what is known about the evolutions of densities in discrete time maps as well as in systems with dynamics described by ordinary differential equations or stochastic differential equations and then highlight the rather formidable mathematical problems that arise when one wishes to consider delayed, or hereditary, dynamics.  \end{quotation}

\section{Introduction}\label{sec:intro}

We are accustomed to thinking about the trajectories of dynamical systems and their possible bifurcations as a parameter is varied.  Typically one encounters bifurcation sequences like
stable steady state $\to$ simple limit cycle $\to$ complicated limit cycle
$\to$ `chaotic' solutions, but the definition of what constitutes chaos is tricky \cite{hunt2015}.

Here, we want to  turn this around and think about the evolution of densities.   This is akin to the Gibbs' notion of looking at an  ensemble of dynamical systems, and this ensemble is described by the corresponding  density  of states. This just means that we are thinking about looking at a very large number of copies of a  dynamical system, under the assumption that each copy is not interacting with any others.

\section{Density evolution in dynamical systems}\label{sec:den}

Though the idea of density evolution may seem an unfamiliar one initially, in point of fact many readers will find that they are really quite familiar with it from other contexts.

From a purely formal standpoint we start with the definition of the Frobenius-Perron (FP) operator $P^t\colon L^1 \to L^1 $
    \begin{equation}
    \int_A P^t f(x)\, \me(dx) =  \int_{\SM_t^{-1}(A)} f(x)
    \,\me(dx)
    \label{eq:FP}
    \end{equation}
which maps densities to densities.  From a technical standpoint\cite{LM94},
$(X,{\cal A},\me)$ is a $\sigma$-finite measure space,  and
$\SM_t\colon X \to X$ a measurable
nonsingular transformation, \ie{} $\SM_t^{-1}(A)\in {\cal A}$ for all $A\in {\cal A}$ and $\me(\SM_t^{-1}(A))=0$ whenever $\me(A)=0$.

This still may look rather unfamiliar, but some examples will smooth the way.  First of all note that if $A = [a,x]$ the Frobenius-Perron  operator becomes
$$
 \int_a^x P^tf(s) \,ds  =    \int_{\SM_t^{-1}([a,x])}f(s) \,ds,
$$
 so
\begin{equation}
 P^tf(x) =
\dfrac{d}{dx} \int_{\SM_t^{-1}([a,x])}f(s) \,ds. \label{eq:1.3}
\end{equation}
For example, with the tent (hat) map
\begin{equation}
    \SM(x) = \left\{
    \begin{array}{ll}
    a x & \qquad \mbox{for}
    \quad x \in
    \left [ 0, \frac 12  \right ) \\
    a (1 - x ) & \qquad \mbox{for} \quad x \in
    \left [ \frac12, 1 \right ],
    \end{array}
    \right.
\label{eq:hat}
\end{equation}
$S_n$ is the $n$th iterate of $S$, $n\in \mathbb{N}$, and
the corresponding Frobenius-Perron operator is the $n$th iterate $P^n$ of the operator
$$
    Pf(x) = \dfrac{1}{a} \left[ f\left( \dfrac{x}{a} \right) + f\left(
    1-\dfrac {x}{a}\right)\right].
$$

In a more familiar vein if we have a system of ordinary differential equations
$$
    \dfrac{dx_i}{dt} = {\cal F}_i(x), \qquad i=1,\ldots ,d,
$$
then from the definition of the Frobenius-Perron operator we can derive the evolution equation for $f(x,t) = P^tf(x)$:
\begin{equation}
    \dfrac {\partial f}{\partial t} = -\sum_{i=1}^d \dfrac {\partial
    (f {\cal F}_i)}{\partial x_i},
\label{eq:gen liou}
\end{equation}
which is just the  generalized Liouville equation\cite{LM94}.

Finally if we have a stochastic differential equation
$$
    {dx}  = {\cal F}(x) dt + \sigma (x) dw(t)
$$
where $w$ is a Wiener process, then the
evolution equation for the density $f(x,t)\equiv P^tf_0(x)$ is the Fokker-Planck equation\cite{gardinerhandbook}
\begin{equation*}
    \frac {\partial f}{\partial t} = -\sum_{i=1}^d \frac {\partial (f{\cal
    F}_i)}{\partial x_i} + \frac {1}{2} \sum _{i,j=1}^d \frac {\partial ^2 (
    a_{ij} f)}{\partial x_i \partial x_j}
    \end{equation*}
where $a_{ij}(x) = \sum_{k=1}^d \sigma_{ik}(x)\sigma_{jk}(x)$.

However, if we have a variable $x$ evolving under the
action of some dynamics described by a differential delay equation
    \begin{equation}
    x'(t)= \epsilon ^{-1} {\cal F}(x(t),x(t-\tau)), \quad
    x(t) = \phi(t) \quad t \in [-\tau, 0],
    \label{eq:0.0}
    \end{equation}
then things are not so clear.  We would like to know how some initial density of the variable $x$ will
evolve in time.

Denote the `density' of $x$ on the interval $[t-\tau,t]$ by $\rho(x,t,t'), \forall t' \in [t-\tau,t]$.
Then we would like to be able to determine an evolution operator ${\cal U}$ such that the
equation
    \begin{equation}
    {\cal U}^t \rho_0(x) = 0, \quad \rho_0(x)  = \rho(x,0,t')    \label{eq:0.1}
    \end{equation}
describes the evolution of $\rho(x,t,t')$ given a density of initial functions $\rho_0(x,0,t')$.
Unfortunately,  we don't really know how to do this, and that's the
whole point of this paper.

The reason that the problem is so
difficult is embodied in \eqref{eq:0.0} and the infinite
dimensional nature of the problem because of the necessity of
specifying the {\it initial function}  $\phi(t)$ for $t \in
[-\tau, 0]$, and even further the `density' of initial functions $\rho_0(x,0,t')$ for $t' \in [-\tau,0]$.

We know about what ${\cal U}^t$   should
look like in various limiting cases.  For instance, to consider an extensively studied example\cite{mallet1986global}, if  ${\cal F}(x(t),x(t-\tau)) = -x(t) + \SM(x(t-\tau))$
so \eqref{eq:0.0} becomes
    \begin{equation}
    \epsilon x'(t)= -x(t) +\SM (x(t-\tau)), \quad
    x(t) = \phi(t) \quad t \in [-\tau, 0],
    \label{eq:0.2}
    \end{equation}
then we expect that:

\begin{enumerate}

\item If we let $\epsilon \to 0$ and restrict consideration to $t/\tau  \in \mathbb{N}$,
then ${\cal U}^t$ as in \eqref{eq:0.1} should reduce to the Frobenius-Perron operator \eqref{eq:1.3}    for the map $\SM $.

\item If $\tau \to 0$ then we should recover the Liouville
equation \eqref{eq:gen liou} from  ${\cal U}^t$.

\item If $\epsilon \to 0$, then from  ${\cal U}^t$ as in \eqref{eq:0.1}  we should recover the
operator governing the evolution of densities {\it in a function
space} under the action of the functional map
    \begin{equation}
    x(t) = \SM (x(t-\tau)),
    \label{eq: 0.3}
    \end{equation}
for $t \in  \mathbb{R}^+$,  though we don't know what that should be\cite{sharkovski-bk}.

We must note that in this case if $t/\tau  \in \mathbb{N}$, then it has been shown\cite{mallet1986bifurcation} that the structure of the dynamics of \eqref{eq: 0.3} is not necessarily preserved by \eqref{eq:0.2}, so what would happen to the questions of `densities' is totally unclear.
\end{enumerate}

\subsection{Classifying density evolution dynamics}\label{ssec:class}

Just as dynamists classify the different types of trajectory behaviours\cite{kuznetsov2013elements}, ergodic theorists have also classified different types of convergence of density evolution (Ref. \onlinecite[Chapter 4]{LM94}).  In this classification, we always let $\SM\colon X \to X$ be a nonsingular transformation on a $\sigma$-finite measure space $(X,{\cal A},\me)$ that preserves
a probability measure $\mu$
with density $f_*$.

The weakest type of convergence is contained in the property of
ergodicity.
$\SM$ is \emph{ergodic} if
every invariant set  $A \in {\cal A}$\footnote{$A$ is invariant if $S^{-1}(A)=A$} is such that either $\mu(A) =
0$ or $\mu(X \setminus A) = 0$.  This is equivalent 
to the existence of a unique stationary density $f_*$ so $Pf_* \equiv f_*$.

Next in the hierarchy is the stronger property of mixing.
$\SM$ is \emph{mixing}  if
    \begin{equation*}
    \lim_{t \to \infty} \mu(A \cap \SM_{t}^{-1}(B)) = \mu(A)\mu(B)
    \quad \mbox {for all } \, A,B \in {\cal A}.
    \label{eq:mixing}
    \end{equation*}
Mixing 
is equivalent to
\begin{equation*}
    \lim_{t\to \infty}   \langle P^tf, g\rangle  = \langle f_*,g\rangle
    \label{eq:mixingdef}
    \end{equation*}
for every bounded measurable function $g$.

Then we have the property of asymptotic stability (or exactness).
Assume $\SM$ is such that $\SM(A) \in {\cal A}$
for each $A \in {\cal A}$. $\SM$ is {\it asymptotically stable} 
if
    \begin{equation*}
    \lim_{t \to \infty} \mu(\SM_t(A) ) = 1 \quad \mbox {for all } \, A,B \in {\cal A}.
    \label{eq:exact}
    \end{equation*}
Asymptotic stability is equivalent to
$$f(x,t)\equiv P^tf_0(x) \to f_*(x)$$
for all initial densities.

Finally there is  asymptotic periodicity (Ref. \onlinecite[Chapter 5.3]{LM94}).  In this case,  for all initial densities $f_0(x)$,   there exists a sequence of basis densities $g_1,\ldots,g_r$ and a sequence of bounded linear functionals $\lambda_1,\ldots,\lambda_r$ such that

$$
Pf_0(x) =  \sum_{i=1}^r \lambda _{i}(f_0)g_{i}(x) + Qf_0(x),
$$
where $Q$ is an operator such that $\parallel P^{t}Q f \parallel \rightarrow 0$ as $t
\rightarrow \infty$ for all integrable $f$.
The densities $g_j$ have  disjoint supports and
$Pg_j=g_{\alpha(j)}$, where $\alpha$ is a permutation of
$(1,\ldots,r)$. An invariant density is given by
$$
    f_*=\dfrac{1}{r}\sum_{j=1}^rg_j.
$$

\subsection*{Example: General hat map}
The hat map is perfect to illustrate these various types of dynamics, since it is known\cite{ito79a,ito79b} that \eqref{eq:hat}
is ergodic for $a > 1$  and we have an analytic expression\cite{yoshida83} for the stationary density $f_*$.  Furthermore, \eqref{eq:hat} is  asymptotically periodic \cite{provatas91a} with period $r = 2^n$,
$n = 0,1,\cdots$ for
$$
    2^{1/2^{ {n+1}} } < a \le 2^{1/2^{ {n}} }.
$$
Thus, for example,  $\{P^tf\}$ has period $1$ for $2^{1/2} < a \leq
2$, period $2$ for $2^{1/4} < a \leq 2^{1/2}$, period $4$ for
$2^{1/8} < a \leq 2^{1/4}$, etc.  Finally it is known \cite{LM94} that \eqref{eq:hat} is exact for $a=2$.

\subsection{Can dynamical systems display a 'chaotic' evolution of densities?}\label{sec:chaotic}

The short answer is that
nobody knows--it's an open problem!

As pointed out in the Introduction, the trajectory sequence of potential solution behaviors  through bifurcations in dynamical or semi-dynamical systems is
\begin{multline*}
\text{stable steady state $\to$  simple limit cycle} \\
\to \text{ complicated limit cycle} \to \text{ `chaotic' solutions}
\end{multline*}
and a great deal is known about the possible transitions between different qualitative behaviours \cite{kuznetsov2013elements}.

Analogously, the bifurcation structure
in the evolution of sequences of  densities under the action of a
Frobenius-Perron operator is\cite{LM94}
\begin{multline*}
\text{asymptotically stable stationary density}\\ \to\text{ simple asymptotic periodicity}\\  \to\text{ complicated asymptotic periodicity}
\end{multline*}
but the analysis of bifurcations of densities is only in a rudimentary state of development (Ref. \onlinecite[Chapter 9]{arnold98}).

Thinking about these two different sequences raises the immediate and obvious question ``How could (can) one construct an evolution
operator for densities that would display a `chaotic' evolution of
densities?''.  Thus, is
\begin{multline*}
\text{asymptotically stable stationary density}\\ \to \text{asymptotic periodicity}\\ \to \text{'chaotic' density evolution}
   \end{multline*}
possible?

The Frobenius-Perron operator \eqref{eq:FP} is a linear operator, so our suspicion is that in order to have a chaotic density evolution
it would be necessary to have a non-linear evolution operator.  We have speculated elsewhere\cite{mackey2009exploring,mackey2012mathematical,mackey2016adventures} that maybe the density evolution operator might need to be density dependent, thus leading to nonlinearity.

\subsection*{Example: A density dependent extension of the hat map}

Consider
a density dependent hat map
    \begin{equation*}
    x_{n+1}=\left\{
    \begin{array}{ll}
    a[f_n] x_n,   & x_n \in [0,\frac 12],\\
    a[f_n](1-x_n), & x_n \in (\frac 12,1 ],
    \end{array}
    \right.
    \end{equation*}
where $f_n$ is the density of $x_n$ and the functional $a[f]$ is defined~by
\begin{equation*}
    a[f] = 1 + \int_A^{A+\delta} f(x) dx.
    \end{equation*}
The corresponding nonlinear evolution (pseudo-Frobenius-Perron)
operator is \citep{mackey2009exploring}
    \begin{equation*}
    P_ff(x) = \dfrac {1_{[0,a[f]/2]}(x)}{a[f]} \left \{ f \left ( \dfrac {x}{a[f]} \right )
    + f \left ( 1 - \dfrac {x}{a[f]} \right )\right \}.
    \end{equation*}

\section{Density evolution in differential delay equations}\label{sec:den-dde}

\subsection{Asymptotic periodicity in a deterministic differential delay equation}\label{ssec:AP deter-delay}

Let $x_{\tau} \equiv x(t-\tau)$ with $\tau = 1$
and consider the hat map \eqref{eq:hat} turned into a delay equation\cite{LM95} (see Eq. \eqref{eq:0.2} in particular)
\begin{equation}
\frac{dx}{dt} = - \alpha  x + \left\{
    \begin{array}{ll}
      ax_\tau & \mbox{if $x_\tau < 1/2$} \\
      a(1-x_\tau) & \mbox{if $x_\tau \geq 1/2$}
 \end{array}
\;\;\;\;\;\;\;\frac {a}{\alpha} \in (1,2], \right.
\label{eq:hat mapp DDE}
\end{equation}
and examine the result of picking many different initial functions and following the trajectories forward in time.  At successive times $t$ we sample across all of the trajectories and form a histogram of the values of $x(t)$ that is an approximation to a `density'.

\begin{figure}
\begin{center}
    \includegraphics[width=0.40\textwidth]{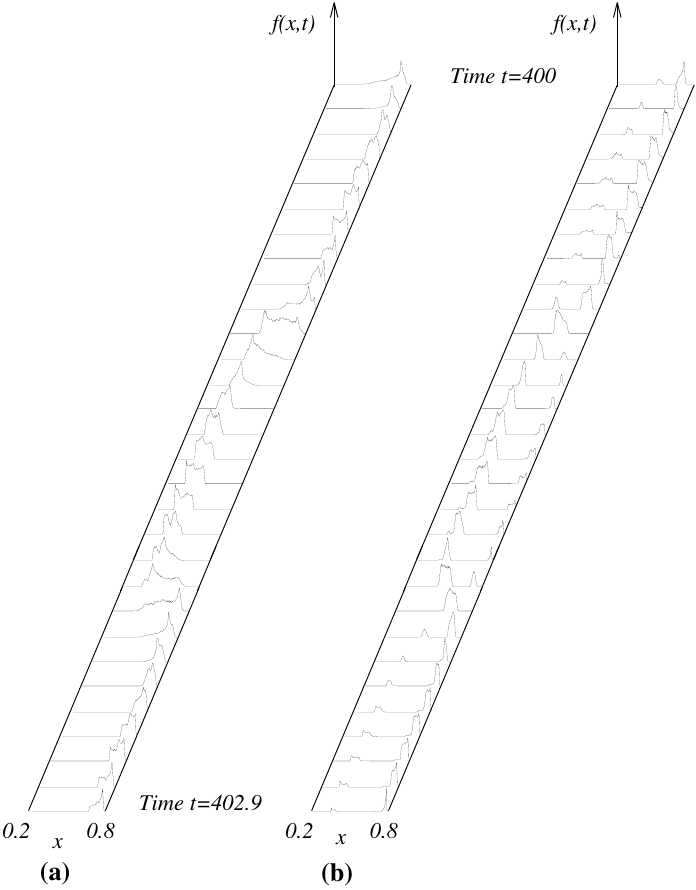}
\end{center}
\caption{Illustration of asymptotic periodicity in the density evolution of a hat map differential delay equation \eqref{eq:hat mapp DDE} with $\alpha = 13$ and $a = 10$.  The simulation extends from $t=400$ to $t = 402.9$ and is based on the integration of $22,500$ initial functions.  (a) Each initial function was a random process distributed uniformly on $[0.65,0.75]$.  (b)  Here the initial functions were uniformly distributed on $[0.65,0.75]$ for $17,000$ of the cases and $[0.35,0.45]$ for the remaining $5,500$.  It was observed that the cycling was not transient and persisted for as long as the simulation ran.  Modified  from Ref. \onlinecite{LM95}.}
 \label{fig:AP dde hat}
\end{figure}

See Figure \ref{fig:AP dde hat} where there is clear numerical evidence for the existence of periodicity in the evolution of the histograms along the trajectories, and which the authors in Ref \onlinecite{LM95} argued was evidence for asymptotic periodicity of densities in this system, supported by their analytic calculations. Note in particular in Figure \ref{fig:AP dde hat} that  a change in the distribution of the initial functions changes the temporal sequence of densities, but not the period.

\subsection{Asymptotic periodicity in stochastically perturbed delay equations}\label{ssec:AP stoch-delay}

\begin{figure}
\begin{center}
    \includegraphics[width=0.50\textwidth]{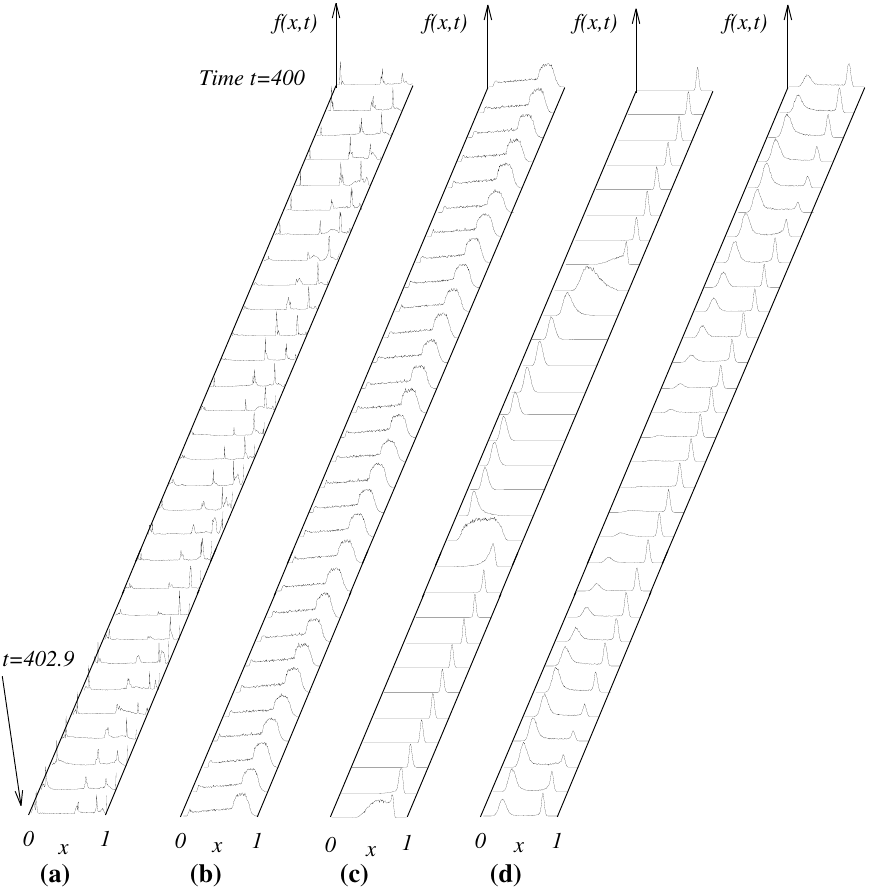}
\end{center}
\caption{Noise induced apparent asymptotic periodicity and stability in the stochastic
  differential delay equation \eqref{eq:keener DDE}.
As in Figure \ref{fig:AP dde hat}, each simulation was performed with
$22,500$ random initial functions. In all four panels, the
parameters of the equation were $a=0.5,\;b=0.567,\;{\alpha}=10$. For
panels (a)-(c) the initial density was as in Figure
\ref{fig:AP dde hat}(a). {(a)} No noise in the system: $f(x,t)$ is
not a density, but rather a generalized function. {  (b)} Noise supported
uniformly on $[0,0.1]$. Numerically, the system appears to be asymptotically stable (exact). {  (c)} Noise uniformly supported
on $[0,0.2]$. {\bf (d)} Same
noise as in (c), with an initial density as in Figure
\ref{fig:AP dde hat}(b).   Modified from Ref. \onlinecite{LM95}.
}
\label{fig:AP dde keener}
\end{figure}

Asymptotic periodicity can be induced by noise\cite{almcm87,provatas91b}  in a Keener map
    \begin{equation*}
    \SM(x) = (ax +b)\,\, \pmod 1, \,\, 0 < a,b < 1,
    \end{equation*}
\ie{} when the dynamics are given by
    \begin{equation}
    x_{n+1} = ( a x_n + b + \xi_n) \,\, \pmod 1, \,\, 0 < a,b < 1,
    \label{eq:noise-keener}
    \end{equation}
and the noise source $\xi$ is distributed with a density $\tilde f$.
Consider the Keener map \eqref{eq:noise-keener} with noise $\xi$ turned into a stochastic delay equation\cite{LM95} (again see Eq. \eqref{eq:0.2})
\begin{equation}
    \frac{dx}{dt} = - \alpha  x + [(ax_\tau +b +\xi) \quad \mbox {mod}\,
    1] \,\, 0<a,b<1 
\label{eq:keener DDE}
\end{equation}
and examine the evolution of many initial functions as shown in Figure \ref{fig:AP dde keener}.

In this figure there are two noteworthy features.  First, an alteration in the distribution of the noise with the distribution of initial functions kept the same leads to an apparent qualitative change in the `density' dynamics, going from asymptotically stable in (b) to asymptotically periodic in (c).  Secondly, a change in the distribution of the initial functions with the distribution of the noise kept constant [(c) to (d)] leads to a change in the details of the temporal density evolution without a change in the period.
\vskip -0.13 in

\subsection{Deterministic Brownian motion}\label{sec:deter BM}
A typical formulation\cite{gardinerhandbook} of the Brownian motion of a particle of mass $m$ with position $x$ and velocity $v$ subject to a frictional force $\gamma v$ is
\begin{align*}
 \dfrac{d x}{d t} &= v, \\ 
m \dfrac{d v}{d t} &= - \gamma v + \eta(t),
\end{align*}
where
$\eta$ is   a fluctuating ``force'' due to collisions of the particle with others of much small mass, and is usually   given by
$
\eta(t) = \sigma \xi(t) $, and $\xi = \frac{dw}{dt}$ is a `white noise' (and delta correlated) which is the `derivative' of a Wiener process $w(t)$.
$\xi(t)$ is normally distributed with mean $\mu = 0$ and variance $\sigma = 1$.

The question of whether one can produce a similar Brownian motion using a totally deterministic model is interesting, and has been answered in the affirmative\cite{lei2011deterministic,mackeytyran}.

\subsection*{Brownian motion from a differential delay equation}\label{ssec:Brown-dde}

\begin{figure}
        \includegraphics[width=0.45\textwidth]{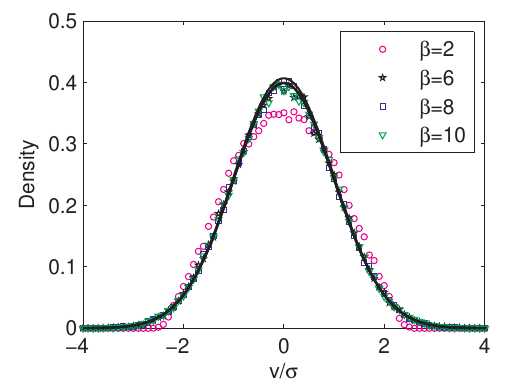}
    \caption{Velocity distribution for four different values of $\beta$ with $\gamma = 1$ when the dynamics are given by \eqref{eq:DDE brown}.  See Ref. \onlinecite{lei2011deterministic} for more details.
    }
    \label{fig:velocity}
\end{figure}

In Ref. \onlinecite{lei2011deterministic}, the authors studied numerically the system
\begin{eqnarray}
\nonumber \dfrac{d x}{d t} &=& v,\\ \label{eq:DDE brown}
 \dfrac{d v}{d t} &=& - \gamma v +   \sin(2 \pi \beta v(t-1)),\\
\nonumber v(t)&=&\phi(t), \,\,-1 \leq t \leq 0.
\end{eqnarray}
In this system the `random' force (the sinusoidal term) is oscillating ever more rapidly as $\beta$ increases.
It was shown \cite{lei2011deterministic}, for a variety of numerical situations, that the mean square displacement   of the particle obeys $[\Delta x(t)]^2 \sim t$, while the velocity is distributed as a quasi-Gaussian $  \sim e^{-Cv^2}$ for $v \in [-K,K]$ (see Figure \ref{fig:velocity}).  The numerics indicated that the bound $K$ and the standard deviation $\sigma$ are given by
\begin{eqnarray*}
\label{eq:K0}
K(\beta,\gamma) &=& \dfrac{1}{\sqrt{\gamma} (0.68\sqrt{\beta} + 0.60 \sqrt{\gamma})},\\
\sigma(\beta, \gamma) &=& \dfrac{0.32}{\sqrt{\beta \gamma}}.
\end{eqnarray*}

\subsection*{Brownian motion induced by perturbation from a non-inver\-tible and chaotic map}\label{ssec:ODE-chaotic-map}

In Ref. \onlinecite{mackeytyran}, the authors examined the properties of the system
\begin{eqnarray*}
 \dfrac{d x}{d t} &=& v,\\
\dfrac{d v}{d t} &=& - \gamma v + \eta(t),
\end{eqnarray*}
where the fluctuating ``force'' consists of a series of delta-function like
impulses given by
\begin{equation*}
\eta(t) = \kappa(\tau)\sum_{j=0}^\infty h(\xi(t)) \delta(t-j\tau)
\end{equation*}
with $\tau>0$ and a scaling parameter $\kappa(\tau)>0$.

The real valued function $h$ is defined on a probability
space $(X,\mathcal{A},\mu)$  and $\xi$ is a highly chaotic deterministic variable defined  by
$\xi(t)=\xi_j $ for $j\tau\le t<(j+1)\tau$, $j\ge 0$,
with $\xi_{j+1}=\SM(\xi_{j})$  and $\xi_0$ the identity on $X$,
where $\SM\colon X\to X$ is an  ergodic map with invariant measure $\mu$, \eg{} $\SM$ is the hat map \eqref{eq:hat} on $[0,1]$ and $h$ is the identity.
It was shown in Ref. \onlinecite[Section 4]{mackeytyran} that   the limit $\tau\to 0$   reproduces the characteristics
of an Ornstein-Uhlenbeck process which is
the solution of the stochastic differential equation
\begin{align*}
dX(t)&=V(t)dt,\\
dV(t)&=-\gamma V(t)dt +\sigma dw(t),
\end{align*}
where $\{w(t):t\ge 0\}$ is a Wiener process, provided
\begin{enumerate}
\item $\kappa(\tau)^2/\tau$ converges to $1$  as $\tau\to
0$ and
\item there is $\beta>1/2$ such that
\[
\limsup_{t\to \infty} t^{2\beta}\|P^t h\|_1<\infty,
\]
where $\|\cdot\|_1$ is the norm in $L^1(X,\mathcal{A},\mu)$ and $P^t$ is the Frobenius-Perron operator \eqref{eq:FP} corresponding to $\SM_t$.
\end{enumerate}
 For precise statements and further examples see \cite{mackeytyran,mtyran14}.

\section{What do these examples show?}\label{sec:what}

\begin{figure*}
\begin{center}
  \includegraphics[width=0.9\textwidth]{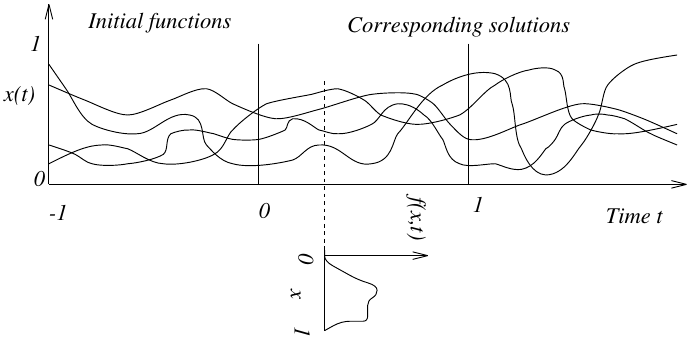}\\
  \end{center}
  \caption{A schematic illustration of the connection between the
  evolution of an ensemble of initial functions and what would be
  measured in a laboratory.  In the case that the delay has been
  scaled to $\tau = 1$, an ensemble of $N$ initial functions on $[-1,0]$ is allowed to evolve forward
  in time under the action of the delayed dynamics. At time $t$ we sample the distribution of
  the values of $x$ across all $N$ trajectories and form an approximation to a density
  $f(x,t)$.
  Modified from Ref. \onlinecite{LM95}.}
  \label{fig:f1d}
\end{figure*}

So what are we displaying in Figures \ref{fig:AP dde hat} through \ref{fig:velocity}?  Reference to Figure \ref{fig:f1d} makes it clear that what we are examining is not really the evolution of the density $\rho(x,t,t')$ but rather $f(x,t) = \rho(x,t,t'\equiv t)$.

This is not necessarily a bad thing because $f(x,t) = \rho(x,t,t'\equiv t)$ is the quantity that would typically be measured in an experimental situation:  We make measurements of $x$ at either discrete times, or a continuum of times, and from these measurements construct temporal histograms of some state variable as in Figures \ref{fig:AP dde hat} or \ref{fig:AP dde keener}, or maybe we look at the long time limiting behaviour as in Figure \ref{fig:velocity}.

In Figure \ref{fig:f1d}
we are looking at a snapshot of the evolution of all of these trajectories emanating from many initial functions just as we would do in an experiment.  However, the unanswered question is how this is related to the {\it density} evolving under the delayed dynamics?  That is, how do we get from $\rho(x,t,t')$ to $f(x,t) = \rho(x,t,t'\equiv t)$?
This is precisely the problem raised in Section \ref{sec:den}.

\section{How to formulate?}\label{sec:formulate}

We now turn to a consideration of the fundamental problem posed in this paper.  Namely, given a system
whose evolution is determined by a differential delay equation of the form \eqref{eq:0.0} for example,
and whose initial phase point $\phi \in C$ is distributed over many possible initial
states, how does the probability distribution for the phase point
evolve in time?  Or equivalently given an  ensemble of independent
systems, each governed by \eqref{eq:0.0},
and whose initial
functions are distributed according to some density over $C$, how does
this ensemble density evolve in time?

The answer to this question {\it seems} deceptively simple and would appear to be provided by
the Frobenius-Perron  operator formalism of Section \ref{sec:den}.
Suppose the initial distribution of phase
points is described by a probability measure $\mu_0$ on $C$ so the probability that the initial function $\phi$ is an element of $A \subset C$
is given by $\mu_0(A)$.
Then, after a time $t$, the new distribution is described
by the measure $\mu_t$ given by
\begin{equation} \label{eq.pullback1}
  \mu_t = \mu_0 \circ \SM_t^{-1}
\end{equation}
if $ \SM_t $ is a measurable transformation on $C$.  Thus after a time $t$ the
probability that the phase point is an element of $A \subset C$ is $\mu_t(A) = \mu_0( \SM_t^{-1}(A))$.

If the initial distribution of states
$\phi$ is described by a density $f(\phi)$ with respect to some
measure $\me$ and $\SM_t$ is a nonsingular transformation, then after time $t$ the density will have evolved
to $P^tf$, where the Frobenius-Perron operator $P^t$ corresponding
to $\SM_t$ is defined by
\begin{equation} \label{eq.fpdde1}
  \int_A P^t f(\phi)\,\me(d\phi) = \int_ {\SM_t^{-1}(A)}
  f(\phi)\,\me(d\phi)
\end{equation}
for all measurable $ A   \subset C.$

Equations \eqref{eq.pullback1}--\eqref{eq.fpdde1} apparently
answer the question posed above about the evolution of probability measures for delay differential equations.
However, {\it this is illusory\cite{losson2020}} since they provide only a symbolic
restatement of the problem.  Everything that is specific to a
given differential delay equation is contained in   $\SM_t^{-1}$.

Although the differential delay equation can be expressed in terms of an evolution semigroup,
there is no apparent way to invert the resulting transformation  $\SM_t$.   This inversion will most certainly be non-trivial, since
solutions of delay equations often cannot be uniquely extended
into the past~\citep{driver2012ordinary}, and thus $\SM_t$ will not have a unique
inverse.  $\SM_t^{-1}$ may have numerous branches that need to
be accounted for when evaluating  $\SM_t^{-1}(A)$ in the Frobenius-Perron
equation~\eqref{eq.fpdde1}.  This is a serious barrier to deriving a
closed-form expression for the Frobenius-Perron operator $P^t$.

There are other subtle issues\cite{losson2020} raised by
equations~\eqref{eq.pullback1}--\eqref{eq.fpdde1}.
\begin{itemize}
\item The most apparent
difficulty is that the integrals in~\eqref{eq.fpdde1} are over sets in
a function space, and it is unclear how such
integrals can be carried out.
\item It is also unclear
what family of measures we are considering, and in particular what
subsets $A \subset C$ are measurable (\ie, what is the relevant
$\sigma$-algebra on $C$?).
\item Also, in equation~\eqref{eq.fpdde1} what
should be considered a natural choice for the measure $\me$ with
respect to which probability densities are to be defined?
\item And finally, does it make sense to talk about probability densities in
the function space $C$?
\end{itemize}

These, then, are the rather formidable mathematical problems that we see in any attempt to formulate a framework to describe the evolution of a `density' under the action of a dynamical systems containing delays.  These have bedeviled us, as well as many of our collaborators, for a number of years.  We now turn to a brief consideration of failed attempts that we have published over the years in an attempt to solve this problem.   A detailed presentation of these can be found in Ref.~\onlinecite{losson2020}.

\subsection*{Possible   lines of attack}
\begin{itemize}
\item Considerations related to an examination of fluid flow, and turbulence in particular, bear some superficial similarities to the problems involved in thinking about density evolution in delay systems because both are infinite dimensional.  In Ref. \onlinecite[Chap. 5]{losson2020} an extension of the work of Ref. \onlinecite{hopf1952statistical} involving Hopf functionals is examined as a possible way of tackling delayed density evolution.  This is a exposition of work first published in Ref. \onlinecite{losson92}, and though a Hopf-like functional differential equation governing the evolution of these densities was derived it was not possible to proceed further with this method.  Interestingly this approach leads to a suggestive potential connection with the functional expansions of quantum field theory and Feynman diagrams.
\item The method of steps\cite{driver2012ordinary}  is a classic tool used in the proof of the existence and uniqueness of solutions of differential delay equations and Ref. \onlinecite[Chap. 6]{losson2020} examined its potential utility for a formulation of the delay differential equation density problem. This approach, while promising and useful from a numerical perspective, seems to lead only to a weak solution and to not be of great utility.

\item As we have noted above, an adequate theory of integration on infinite
dimensional spaces is lacking.  Such a theory is needed if we are to
further develop the Frobenius-Perron operator formalism to
characterize the evolution of densities for delay differential equations, which requires a
theory of integration of functionals.
However, there is a notable exception worth mentioning.  There is one
probability measure (or family of measures) on a function space,
called Wiener measure, for which there is a substantial theory of
integration\cite{Kac}.

With Wiener measure it has been possible to prove strong ergodic properties
(\eg~exactness) for a certain class of partial differential
equations\cite{BK84,Rud85,Rud87,Rud88}.  The success of these
investigations, together with the considerable machinery that has been
developed around the Wiener measure, suggests that Wiener measure
might be a good choice for the measure of integration in the study of
other infinite dimensional systems such as delay equations but this requires further investigation.

\item Another possible approach is to take a differential delay equation and, through a process of discretization (e.g. with an Euler approximation) turn it into a high dimensional map\cite{losson1995evolution}.  Using this approach it was then possible to prove certain properties for the density evolution of the high dimensional map using results from Ref. \onlinecite{ionescu50}, but it was never possible to successfully make the transition from these results to a consideration of the discretized differential delay equation (Ref. \onlinecite[Chaps. 7,8]{losson2020}).

\item In the same spirit of the previous suggestion, another numerical  approximation approach\cite{dellnitz1999approximation} has been extended\cite{dellnitz2015computation} to delay equations to obtain a representation of the limiting attractor in a projected space.  Numerical efforts like these are extremely valuable because they offer potential insights into what analytic approaches should be giving as answers once we figure out what the analytic approach actually is.  Many more examples of the utility and the shortcomings of this numerical approximation approach can be found in Ref. \onlinecite{losson2020}.

\end{itemize}

None of these approaches have yielded a satisfactory way of dealing with the problems enumerated at the beginning of this section and it is clear from the extensive considerations of Ref. \onlinecite{losson2020} that a radically new approach is required.  In the next section we outline a possible approach and show its application to a very simple example.

\section{A novel approach}\label{s:weak}

In this section we outline a new approach to the delayed dynamics density evolution problem formulated in \eqref{e:dde} and show its potential applicability to a specific and relatively simple, tractable, example.

Briefly in the first part of this section ({\it Liouville-like formulation}) we straightforwardly derive an evolution equation \eqref{e:weakL} for $f(x,t)$ when the dynamics are given by \eqref{e:dde}, and note that the equation that we obtain reduces to the standard Liouville equation \eqref{eq:gen liou} when $\mathcal{F}$ does not depend on $x(t-\tau)$.

We then turn in the second part of the section (entitled {\it A tractable example:  Gaussian processes}) to a consideration of a Gaussian process $\xi$ with continuous sample paths on the initial interval $[-\tau,0]$ as the initial function for the simple linear differential delay equation \eqref{e:lin}.  We derive the analog of \eqref{e:weakL} in \eqref{e:weakL1}, and are able to analytically characterize the nature of the solutions for all $(a,b,\tau)$.

\subsection*{Liouville-like formulation}

Let $C=C([-\tau,0],\mathbb{R}^d)$ denote the Banach space of all continuous
functions $\phi\colon[-\tau,0]\to\mathbb{R}^d$ equipped with the supremum norm and the Borel $\sigma$-algebra.
Consider the equation
\begin{equation}\label{e:dde}
\begin{split}
   x'(t)=&\mathcal{F}(x(t),x(t-\tau)), \quad t\ge 0, \\
    x(t)=&\phi(t),\quad t\in[-\tau,0],
  \end{split}
\end{equation}
where $\phi\in C$ and $\mathcal{F}$ is such that for each $\phi\in C$  there exists a continuous function
$x\colon[-\tau,\infty)\to \mathbb{R}^d$ such that \eqref{e:dde} has a unique
global solution depending continuously on~$\phi$. For each $t\ge 0$ define the solution map $\SM_t\colon
C\to C$ by
\begin{equation}\label{d:sol}
\SM_t(\phi)(s)=x_t(s)=x(t+s)\quad \text{for }
s\in[-\tau,0],
\end{equation}
where $x(t)$ is a solution of \eqref{e:dde} with $x_0=\phi$. It is
well known\cite{hale-lunel} that $\{\SM_t\}_{t\in\mathbb{R}^+}$ is a semi-dynamical system on
$C$ and the transformation $(t,\phi)\mapsto \SM_t(\phi)$ is
continuous.

Let $B(C)$ be the space of bounded Borel measurable functions $\psi \colon C\to \mathbb{R}$
with the supremum norm
\[
\|\psi \|_{\infty}=\sup_{\phi\in C}|\psi (\phi)|,
\]
and let $\mathcal{M}(C)$ (resp. $\mathcal{M}_1(C)$) denote the space of finite (resp. probability) Borel measures on $C$.
For any $\psi \in B(C)$ and $\mu\in \mathcal{M}(C)$ we use the scalar product notation
\[
\langle \psi ,\mu\rangle=\int_C \psi (\phi)\mu(d\phi).
\]
Let $C_b\subset B(C)$ be the Banach space of all bounded uniformly continuous
functions $\psi\colon C\to \mathbb{R}$ with the supremum norm. We say that a family
$\psi_t\in C_b$, $t>0$ \emph{converges weakly}\cite{dynkin,mohammed} to $\psi\in C_b$  as $t\to 0^{+}$ (denoted by
$\mathrm{w-}\lim_{t\to 0}\psi_t=\psi$) if
\[
\lim_{t\to 0}\langle \psi_t,\mu\rangle=\langle \psi,\mu\rangle\quad \text{for each }
\mu\in \mathcal{M}(C).
\]
This is equivalent to the following two conditions: $\lim_{t\to
0}\psi_t(\phi)=\psi(\phi)$ for every $\phi\in C$  and
$\sup_t\norm{\psi_t}_{\infty}<\infty$.

A semigroup $\{T^t\}_{t\ge
0}$ of linear operators on the space $C_b$ is defined by
\begin{equation}\label{eq:sem}
T^t\psi(\phi)=\psi(\SM_t(\phi)),\quad t\ge 0,\;\phi\in C, \psi\in
C_b,
\end{equation}
where $\SM_t$ is the solution map \eqref{d:sol}.
 Since the semidynamical system $\{\SM_t\}_{t\ge 0}$ is
continuous, we obtain
\[
\lim_{t\to 0}T^t\psi(\phi)=\psi(\phi)\quad \text{for all } \phi\in C,\; \psi\in
C_b.
\]
Note also that
\[
\sup_{t\ge 0}\|T^t\psi\|_{\infty}\le \|\psi\|_{\infty}.
\]
Consequently, the semigroup $\{T^t\}_{t\ge 0}$ is weakly continuous at $t=0$.
Define the \emph{weak generator} $\mathcal{L}\colon \mathcal{D}(\mathcal{L})\subset
C_b\to C_b$ of the semigroup $\{T^t\}_{t\ge 0}$ by\cite{dynkin,mohammed}
\begin{align*}
\mathcal{D}(\mathcal{L})&=\{\psi\in C_b: \mathrm{w-}\lim_{t\to
0}\frac{1}{t}\bigl(T^t\psi-\psi\bigr)\; \text{exists}\},\\
\mathcal{L}\psi &=\mathrm{w-}\lim_{t\to 0}\frac{1}{t}\bigl(T^t\psi-\psi\bigr).
\end{align*}
In particular, for $\psi \in \mathcal{D}(\mathcal{L})$ and $\mu \in \mathcal{M}(C)$ we have
\begin{equation}\label{e:weak}
\langle T^t \psi,\mu\rangle=\langle  \psi,\mu\rangle +\int_0^t \langle \mathcal{L} (T^r\psi),\mu\rangle  dr,\quad t>0.
\end{equation}

Let $\mu_0\in \mathcal{M}_1(C)$. For each $t\ge 0$, define  the probability measure $\mu_t$ on the space $C$ by $\mu_t=\mu_0\circ \SM_t^{-1}$ as in  \eqref{eq.pullback1},
where $\SM_t$ is the solution map \eqref{d:sol}. Then, by \eqref{eq:sem}, we have
\[
\langle \psi,\mu_t\rangle =\langle T^t \psi,\mu_0\rangle,\quad t\ge 0, \psi\in C_b.
\]
It follows from \eqref{e:weak} and the change of variables formula that
\begin{equation}\label{e:weakf}
\langle \psi,\mu_t\rangle=\langle  \psi,\mu_0\rangle +\int_0^t \langle \mathcal{L} \psi,\mu_r\rangle  dr
\end{equation}
for all $ t>0$ and $ \psi \in \mathcal{D}(\mathcal{L})$. However, it is difficult to identify the domain $ \mathcal{D}(\mathcal{L})$ of the weak generator $ \mathcal{L}$.
Thus, we introduce an \emph{extended generator} for the solution map $\SM_t$. It is defined as  a linear operator $\mathcal{L}$
from its domain $\mathcal{D}$ to the set of all Borel measurable functions on $C$, where we say that $\psi \in \mathcal{D}$ if  for each $t>0$ we have
\begin{equation*}
\int_0^t \langle |\mathcal{L} \psi|,\mu_r\rangle  dr <\infty
\end{equation*}
and \eqref{e:weakf} holds. 
Then we say that $\{\mu_t\}_{t\ge 0}$ is the solution of the equation
\begin{equation}\label{eq:infdim}
\frac{\partial}{\partial t}\mu_t=\mathcal{L}^*\mu_t.
\end{equation}

Instead of defining the whole domain $\mathcal{D}$, one can start with a sufficiently large subset of $\mathcal{D}$. A typical example is $\mathcal{F}C_c^\infty(C)$ the set of smooth cylinder functions on $C$. These are functions of the form
$\psi(\phi)=g(l_1(\phi),\ldots,l_n(\phi))$, $\phi\in C$, $n\ge 1$, for some continuous linear functionals $l_1,\ldots,l_n$ on $C$, and $g\in C_c^\infty(\mathbb{R}^n)$. Another example of $\mathcal{D}$ is the set of quasi-tame functions introduced by Mohammed~\cite{mohammed}. We will use below yet another subset of~$\mathcal{D}$ that will allow us to change \eqref{eq:infdim} into a partial differential equation \eqref{e:weakL}.

The marginal distribution $\mu(t)$ of the measure $\mu_t=\mu_0\circ \SM_t^{-1}$ is defined on $\mathbb{R}^d$  by
\[
\mu(t)(B)=\mu_t\{\phi\in C: \phi(0)\in B\}, \quad B\in \mathcal{B}(\mathbb{R}^d).
\]
This can be rewritten with the projection map $\pi_0\colon C\to \mathbb{R}^d$ defined by $\pi_0(\phi)=\phi(0)$, $\phi\in C$, as
 $\mu(t)=\mu_t\circ \pi_0^{-1}$. Note that $\mu(t)$ is the distribution of $x(t)$ for all $t\ge 0$.

Let $C_c^1(\mathbb{R}^d)$ denote the space of functions that are continuously differentiable and have compact support.
Consider the differential operator from the space $C_c^1(\mathbb{R}^d)$ to the set of all Borel measurable functions on $C$  defined by
\[
Lg(\phi)=\sum_{i=1}^d \mathcal{F}_i(\phi(0),\phi(-\tau))\frac{\partial g}{\partial x_i} (\phi(0)),\quad \phi\in C,g\in C_c^1(\mathbb{R}^d).
\]
Let $g\in C_c^1(\mathbb{R}^d)$ be such that $Lg\in B(C)$. Now if $\psi=g\circ \pi_0$, then $\psi \in \mathcal{D}$ and $\mathcal{L}\psi=Lg$.

We say that $\{\mu_t\}_{t\ge 0}$ is a solution of the equation
\begin{equation*}
\frac{\partial}{\partial t}\mu(t)=L^*\mu_t
\end{equation*}
if for each $ t> 0$ and $ g\in C_c^1(\mathbb{R}^d)$ the following holds
\begin{equation*}
\int_{0}^t \langle |Lg|,\mu_r\rangle dr<\infty
\end{equation*}
and
\begin{equation}\label{e:gLd}
\int_{\mathbb{R}^d} g(x)\mu(t)(dx)=\int_{\mathbb{R}^d} g(x)\mu(0)(dx)+\int_{0}^t \langle Lg,\mu_r\rangle dr.
\end{equation}
Note that this is equivalent to requiring that $\psi\in \mathcal{D}$ and $\mathcal{L}\psi=Lg$ for $\psi=g\circ \pi_0$ and all $g\in C^1_c(\mathbb{R}^d)$. Next, observe that if we introduce the measure
\[
\nu(t)=\mu_t\circ \pi_{0,-\tau}^{-1},\quad t\ge 0,\] where $\pi_{0,-\tau}\colon C\to \mathbb{R}^d\times \mathbb{R}^d$ is the projection map $\pi_{0,-\tau}(\phi)=(\phi(0),\phi(-\tau))$, $\phi\in C$, then
\[
\langle Lg,\mu_t\rangle=\int_{\mathbb{R}^d}\int_{\mathbb{R}^d}\sum_{i=1}^d \mathcal{F}_i(x,y)\frac{\partial g}{\partial x_i} (x)\nu(t)(dx,dy),
\]
by the change of variables formula. The measure $\nu(t)$ is the distribution of $(x(t),x(t-\tau))$.

Now suppose additionally that the measure $\nu(t)$ has a density $f_\nu$ with respect to the Lebesgue measure on $\mathbb{R}^{2d}$ \ie, $\nu(t)(dx,dy)=f_\nu(x, y,t)dxdy$. Then the measure $\mu(t)$ has a density $f$ with respect to the Lebesgue measure on $\mathbb{R}^d$ and
\[
f(x,t)=\int_{\mathbb{R}^d}f_\nu(x,y,t)dy.
\]
We also have
\[
f(y,t-\tau)=\int_{\mathbb{R}^d}f_\nu(x,y,t)dx, \quad t\ge \tau.
\]
We can rewrite
\eqref{e:gLd} as
\begin{align*}
&\int_{\mathbb{R}^d} g(x)f(x,t)dx=\int_{\mathbb{R}^d} g(x)f_0(x)dx\\
&+\int_0^t \int_{\mathbb{R}^d} \sum_{i=1}^d \mathcal{F}_i(x,y)\frac{\partial g}{\partial x_i} (x)f_\nu(x,y,r)dxdy dr,
\end{align*}
which is the weak form of the equation
\begin{equation}\label{e:weakL}
\frac{\partial}{\partial t} f(x,t)=-\sum_{i=1}^d \frac{\partial}{\partial x_i}\int_{\mathbb{R}^d} \mathcal{F}_i(x,y)f_\nu(x,y,t)dy.
\end{equation}
Eq. \eqref{e:weakL} reduces to the  Liouville equation if $\mathcal{F}$ does not depend on $y$.

Finally, recall that a Borel probability measure $\mu$ on $C$ is called a Gaussian measure (see Ref. \onlinecite{bogachev98}) if the measure $\mu\circ l^{-1}$ is Gaussian on $\mathbb{R}$ for each continuous linear functional~$l$ on~$C$.
Suppose that the solution map $S_t$ is linear. Then if we take as $\mu_0$ a Gaussian measure on $C$, the measure $\mu_t$ will be again Gaussian.
We will look at such examples next.

\subsection*{A tractable example: Gaussian processes}
A Gaussian process is a family $\xi=\{\xi(t)\}_{t\in \mathbb{T}}$ of (real-valued) random variables defined on some probability space $(\Omega,\Sigma,\mathbb{P})$, indexed by a parameter set $\mathbb{T}$, such that every finite linear combination $\sum c_{t_i}\xi(t_i)$ is either identically zero or has a Gaussian distribution on $\mathbb{R}$. Given a Gaussian process  $\xi$ its mean function is  $a(t)=\mathbb{E}(\xi(t))$, $t\in \mathbb{T}$, and its covariance function is the bivariate symmetric function
\[
R(t,s)=\mathrm{cov}(\xi(t),\xi(s))=\mathbb{E}(\xi(t)-a(t))(\xi(s)-a(s)).
\]
The process is called centered if its mean function is zero. Note that a covariance function is  non-negative definite, i.e.
\begin{equation}\label{d:nnd}
\sum_{i,j=1}^n R(t_i,t_j)c_{i}c_{j}\ge 0
\end{equation}
for any $t_1,\ldots,t_n\in \mathbb{T}$ and $c_{1},\ldots,c_{n}\in\mathbb{R}$, $n\ge 1$.  A Gaussian process is said to be non-degenerate if its covariance function is  positive definite, i.e.  inequality \eqref{d:nnd} is strict  for all non-zero $(c_{1},\ldots,c_{n})\in \mathbb{R}^n$.

For a centered Gaussian process the Gaussian distribution of a finite   combination $\sum c_{i}\xi(t_i)$ is  determined through the variance $\mathbb{E}(\sum c_{i}\xi(t_i))^2$ that can be calculated with the help of the covariance function $R$.  Thus the
covariance function of
a centered Gaussian process completely determines all of the finite-dimensional distributions (that is, the
joint distributions of $(\xi(t_1),\ldots,\xi(t_n))$ for any $t_1,\ldots,t_n\in \mathbb{T}$ and $n\ge 1$).
Consequently the distribution of the entire centered Gaussian process $\xi$ is uniquely determined through its covariance function.

Given a symmetric non-negative definite function $R\colon [-\tau,0]\times[-\tau,0]\to\mathbb{R}$ there exists a 
centered Gaussian process with $R$ being its covariance function, by the Kolmogorov extension theorem. If the process has continuous sample paths then it defines a Gaussian measure on~$C$, the distribution of the process $\xi$, by
\[
\mu(B)=\mathbb{P}(\xi\in B), \quad B\in \mathcal{B}(C).
\]
If we have two centered Gaussian processes with continuous sample paths and the same covariance function then they have the same distribution.

Consider the Gaussian process
\[
\xi(t)=\zeta \cos(t-\theta), \quad t\in \mathbb{R},
\]
where $\zeta$ and $\theta$ are independent random variables, the amplitude $\zeta$ has the Rayleigh distribution with density $x e^{-x^2/2}$, $x\ge 0$, and $\theta$ is uniformly distributed on $[0,2\pi)$. Observe that the mean  $\mathbb{E}(\xi(t))$ is zero and the covariance function is
\[
R(t,s)=\mathrm{cov}(\xi(t),\xi(s))=\mathbb{E}(\xi(t)\xi(s))=\cos(t-s).
\]
Then $\xi(t)$ is normally distributed with mean $0$ and variance $1$. Since $R(t,t-\tau)=0$ with $\tau=\pi/2$ and $(\xi(t),\xi(t-\tau))$ has a Gaussian distribution on $\mathbb{R}^2$, we see that $\xi(t)$ and $\xi(t-\tau)$ are uncorrelated, thus independent. Note that $\xi$ is a solution of
\begin{equation}\label{e:cos}
x'(t)  = -x\Big(t - \frac{\pi}{2}\Big).
\end{equation}
We take as $\mu_0$  the Gaussian measure being the distribution of the Gaussian process $\xi=\{\xi(s)\}_{s\in [-\tau,0]}$.  Then $S_t\xi$ is a  centered Gaussian process with the same covariance function as $\xi$. Hence $\mu_t=\mu_0\circ S_t^{-1}=\mu_0$ for all $t\ge0$ implying that the measure $\mu_0$ is invariant for the solution map $S_t$ of equation \eqref{e:cos}.  We also have $f(x,t)=f(x)$ and   $f_\nu(x,y,t)=f(x)f(y)$, where $f(x)$ is the density of the standard normal distribution.   Hence  \eqref{e:weakL} holds with $d=1$ and $\mathcal{F}(x,y)=-y$. However, if we take $\xi(t)=\zeta \cos(t)$, $t\in \mathbb{R}$, where now $\zeta$ has the standard normal distribution, then $\xi$ is a Gaussian process with $R(t,s)=\cos(t)\cos(s)$ and $\xi(t)$ is also a solution of \eqref{e:cos}, but $\xi$ is degenerate since the distribution of  $(\xi(t),\xi(t-\tau))=\zeta(\cos(t),\sin(t))$ is concentrated on a line.

We can extend  \eqref{e:cos} in the following way.
Consider the linear differential delay equation
\begin{equation}\label{e:lin}
x'(t)=ax(t)+bx(t-\tau),
\end{equation}
where $a,b$ are real constants
and  a Gaussian process $\xi=\{\xi(s)\}_{s\in [-\tau,0]}$ with continuous sample paths on $[-\tau,0]$ as the initial condition $x_0=\xi$.
Then the solution map $S_t$ is a linear mapping of $C$, and thus  $S_t\xi$ is a Gaussian process, so  the distribution of $S_t\xi$ is a Gaussian measure on $C$.

Thus we can start with a positive definite symmetric function $R_0$ such that the corresponding Gaussian process $\xi$ has continuous sample paths on $[-\tau,0]$. Then $S_t\xi$ will be a centered Gaussian process with a covariance function $R_t$.
The measure $\mu(t)$, being the distribution  of $x(t)=S_t\xi(0)$, is Gaussian with mean zero $0$ and variance
$\sigma^2(t)=R_t(0,0)$
and the measure $\nu(t)$, being  the distribution of $(x(t),x(t-\tau))=(S_t\xi(0),S_t\xi(-\tau))$, is Gaussian with covariance matrix
\[
Q_t=\left(
  \begin{array}{cc}
    R_t(0,0) & R_t(-\tau,0) \\
    R_t(-\tau,0) & R_t(-\tau,-\tau) \\
  \end{array}
\right).
\]
Note that $\sigma^2(t-\tau)=R_0(t-\tau,t-\tau)$ for $t\in [0,\tau)$ and
 $\sigma^2(t-\tau)=R_t(-\tau,-\tau)$ for $t\ge \tau$.
Now if $\sigma^2(t)\neq 0$ then the density $f$ of $\mu(t)$ is
\begin{equation}\label{e:densf}
f(x,t)=\frac{1}{\sqrt{2\pi\sigma^2(t)}}e^{-\frac{x^2}{2\sigma^2(t)}}
\end{equation}
and if $\det (Q_t)\neq 0$ then the density $f_\nu$ of $\nu(t)$ is given by
\begin{equation}\label{e:densfn}
f_v(x,y,t)=\frac{1}{2\pi \sqrt{\det(Q_t)}}e^{-\frac{\sigma^2(t-\tau)x^2+2R_t(-\tau,0)xy+\sigma^2(t)y^2}{2\det(Q_t)}}.
\end{equation}
Observe that in this example Eq. \eqref{e:weakL} is of the form
\begin{equation*}
\frac{\partial}{\partial t} f(x,t)=-\frac{\partial}{\partial x}(ax f(x,t))-\frac{\partial}{\partial x}\int_{\mathbb{R}} byf_\nu(x,y,t)dy
\end{equation*}
and that
\[
\int_{\mathbb{R}} yf_\nu(x,y,t)dy=\frac{R_t(-\tau,0)}{\sigma^2(t)}xf(x,t).
\]
Thus $f$ is a solution of
\begin{equation}\label{e:weakL1}
\frac{\partial}{\partial t} f(x,t)=-\Big(a+b\frac{R_t(-\tau,0)}{\sigma^2(t)}\Big)\frac{\partial}{\partial x}(xf(x,t)).
\end{equation}
It is easily seen that $f$ satisfies Eq. \eqref{e:weakL1} if and only if
\begin{equation}\label{e:sigr}
\frac{d}{dt}(\sigma^2(t))=2a\sigma^2(t)+2bR_t(-\tau,0),\quad t>0.
\end{equation}
It follows from \eqref{e:lin} that \eqref{e:sigr} holds.
We see that
\[
\sigma^2(t)=e^{2at}\sigma^2(0) +2b \int_0^t e^{2a(t-s)} R_s(-\tau,0)ds,\quad t>0.
\]
Consequently, to obtain the densities $f$ and $f_{\nu}$ it is enough to determine first $R_t(-\tau,0)$ and then $\sigma^2(t)$.

To find the covariance function $R_t$ observe that we can write the solution of \eqref{e:lin} as (see Section 1.6 in Ref. \onlinecite{hale-lunel})
\[
x(t)=X(t)\xi(0)+b\int_{-\tau}^0 X(t-r-\tau)\xi(r)dr,
\]
where $X(t)$ is the fundamental solution of  \eqref{e:lin}, i.e.  $X(t)=0$ for $t<0$ and
\begin{equation}\label{e:funs}
X(t)=\sum_{k=0}^{\lfloor t/\tau\rfloor} \frac{b^k}{k!}(t-k\tau)^k e^{a(t-k\tau)},\quad t\ge 0,
\end{equation}
where $\lfloor s\rfloor=\max\{k\in \mathbb{Z}: k\le s\}$.
We rewrite $x(t)$  with the help of the Lebesgue-Stieltjes integral
as
\begin{equation*}
x(t)=\int_{-\tau}^0 \xi(r) dM_t(r),
\end{equation*}
where the function $M_t\colon [-\tau,0]\to \mathbb{R}$ is defined by
\[
dM_t(r)=bX(t-r-\tau)dr+X(t)\delta_{0}(dr),\quad t\ge 0,
\]
and
\[
dM_t(r)=\delta_{t}(dr),\quad t<0,
\]
with $\delta_{t}$ denoting the point measure at $t$.
Since
\begin{equation}\label{e:St}
S_t \xi(s)=\int_{-\tau}^0 \xi(r) dM_{t+s}(r),
\end{equation}
we obtain
\[
S_t \xi(s_1)S_t \xi(s_2)=\int_{-\tau}^0\int_{-\tau}^0 \xi(r_1)\xi(r_2) dM_{t+s_1}(r_1)dM_{t+s_2}(r_2).
\]
Taking the expectation on both sides of the above leads to
\begin{equation}\label{e:Rt}
R_t(s_1,s_2)=\int_{-\tau}^0\int_{-\tau}^0 R_0(r_1,r_2) dM_{t+s_1}(r_1)dM_{t+s_2}(r_2).
\end{equation}
Since the covariance function $R_t\colon [-\tau,0]\times[-\tau,0]\to \mathbb{R}$ is  symmetric, we can assume that $-\tau\le s_1\le s_2\le 0$.
If $t\in [0,-s_2]$ then
\[
R_t(s_1,s_2)=R_0(t+s_1,t+s_2),
\]
for $t\in (-s_2, -s_1]$ we have
\begin{align*}
R_t(s_1,s_2)&=e^{a(t+s_2)}R_0(0,t+s_1)\\
&\quad +b\int_{-\tau}^{t+s_2-\tau}e^{a(t+s_2-r)}R_0(r, t+s_1)dr
\end{align*}
and if
$t> -s_1$ then
\begin{align*}
R_t(s_1,s_2)&=X(t+s_1)X(t+s_2)R_0(0,0)\\&\quad +bX(t+s_1)\int_{-\tau}^0X(t+s_2-r-\tau)R_0(r,0)dr\\
&\quad +bX(t+s_2)\int_{-\tau}^0 X(t+s_1-r-\tau)R_0(r,0)dr\\
&\quad + b^2\int_{-\tau}^0\int_{-\tau}^0 X(t+s_1-r_1-\tau)X(t+s_2-r_2-\tau)\\
&\quad\quad\quad  \times R_0(r_1,r_2)dr_1dr_2.
\end{align*}
Thus we can find the covariance function $R_t$ by specifying the covariance function $R_0$ and using \eqref{e:funs} in the above equation.

Define
\[
\alpha_0=\max\{\mathrm{Re}\lambda: \lambda=a +be^{-\lambda \tau}\}.
\]
By Theorem 5.2 in Chapter 1 of Ref. \onlinecite{hale-lunel} for each $\alpha>\alpha_0$ there is a constant $c$ such that $|X(t)|\le c e^{\alpha t}$ for all $t>0$. In particular, if $\alpha_0<0$ then $X(t)$ converges to zero exponentially fast.
Since $R_0$ being a continuous function is bounded, we see that $R_t$ approaches $0$ if $\alpha_0<0$.
Based on the work of \cite{hayes1950roots}, see also Ref.  \onlinecite[Section 5.2 and Thm. A.5]{hale-lunel}, we have $\alpha_0<0$
if and only if
\begin{align*}
a\tau <1,\quad  b\tau + a\tau<0,\quad  b\tau+a\tau \cos\kappa +\kappa\sin \kappa>0,
\end{align*}
where $\kappa$ is the root of $\kappa=a\tau \tan \kappa$, $0<\kappa<\pi$ if $a\neq 0$ and $\kappa=\pi/2$ if $a=0$.
These are values  of $(-a,-b)$ inside the cusp like area of Ref. \onlinecite[Figure 5.1]{hale-lunel}. Then we will have
\[
\lim_{t\to\infty} \sigma^2(t)=0
\]
leading to
\[
\lim_{t\to\infty}\int_{-\varepsilon}^\varepsilon f(x,t)dx=1\quad \text{for all }\varepsilon>0,
\]
by Chebyshev's inequality
\[
1-\int_{-\varepsilon}^\varepsilon f(x,t)dx=\mathbb{P}(|x(t)|>\varepsilon)\le \frac{1}{\varepsilon^2}\sigma^2(t).
\]
If, for example $R_0$ is nonnegative with $\sigma^2(0)=R_0(0,0)>0$ and $a\ge 0$, $b>0$ then $\sigma^2(t)\ge X(t)^2 \sigma^2(0)$ and  $X(t)\to \infty$ as $t\to\infty$. Thus
\[
\lim_{t\to\infty} \sigma^2(t)=\infty
\]
implying that
\[
\lim_{t\to\infty}\int_{-\varepsilon}^{\varepsilon}f(x,t)dx=0,
\]
since
\[
f(x,t) \le \frac{1}{\sqrt{2\pi \sigma^2(t)}},\quad x\in \mathbb{R}, t>0.
\]
It might also happen that  $\sigma^2(t)$ is a constant, as it was for Eq. \eqref{e:cos} and $R_0(s_1,s_2)=\cos(s_2-s_1)$.

Suppose now that the covariance function $R_0$ can be written in the form
\begin{equation}\label{e:R0}
R_0(s_1,s_2)=\int_{-\tau}^{0}\eta_r(s_1)\eta_r(s_2)dr,
\end{equation}
for some  function $\eta_r\colon [-\tau,0]\to \mathbb{R}$ such that $(s,r)\mapsto \eta_r(s)$ is Borel measurable.
Then it follows  from \eqref{e:Rt} and \eqref{e:St} that the covariance $R_t$ is given by
\begin{equation}\label{e:Rt1}
R_t(s_1,s_2)=\int_{-\tau}^0 S_t\eta_r(s_1)S_t \eta_r(s_2)dr
\end{equation}
and in particular, we have
\[
\sigma^2(t)=\int_{-\tau}^0 (S_t\eta_r(0))^2dr.
\]
One example of \eqref{e:R0} is
$R_0(s_1,s_2)=\min\{s_1,s_2\}+\tau$ for $s_1,s_2\in [-\tau,0]$, since  \eqref{e:R0} holds with
\begin{equation}\label{e:qrW}
\eta_r(s)=1_{[r,0]}(s),\quad s,r\in [-\tau,0].
\end{equation}
Then we have  $\xi(s)=W(s+\tau)$, where $W=\{W(t)\}_{t\ge 0}$ is the standard Wiener process on $[0,\infty)$,
Another one is
\[
R_0(s_1,s_2)=\left\{
               \begin{array}{ll}
                 u(s_1)v(s_2), &  s_1\le s_2, \\
                 u(s_2)v(s_1), & s_1>s_2,
               \end{array}
             \right.
\]
where $u,v$ are nonnegative differentiable functions with $u(-\tau)=0$ and $v(s)>0$. Here we can define $\xi(s)=v(s) W(\frac{u(s)}{v(s)})$ and we get \eqref{e:Rt1} with
\[
\eta_r(s)=v(s)\sqrt{\frac{u(r)}{v(r)}}1_{[-\tau,s]}(r)\quad s,r\in [-\tau,0].
\]

Finally, we calculate the covariance $R_t$ as in \eqref{e:Rt1} for  $t\in (0,\tau]$, when $\eta_r$ is given by \eqref{e:qrW}, so that
$S_t\eta_r(s)=1_{[r,0]}(t+s)$ for $t+s<0$ and
\[
S_t\eta_r(s)=X(t+s)+b\int_{r}^0 X(t+s-q-\tau)dq
\]
for $s,r\in [-\tau,0]$, $t>-s$.
Observe that
\begin{align*}
S_t\eta_r(s)=\left\{
         \begin{array}{ll}
          1+b(t+s-r-\tau)_{+}, & a=0, \\
 e^{a(t+s)}+\frac{b}{a}(e^{a(t+s-r-\tau)_+}-1), & a\neq 0,
         \end{array}
       \right.
\end{align*}
for  $t+s\in [0,\tau]$, where $(q)_+=\max\{0,q\}$.
Let $-\tau\le s_1\le s_2\le 0$.
We have
\[
R_t(s_1,s_2)=t+s_1+\tau, \quad t\in[0,-s_2],
\]
but if $t\in (-s_2, -s_1]$  then
\[
R_t(s_1,s_2)=\left\{
         \begin{array}{ll}
          t+s_1+\tau +\frac{b}{2}(t+s_2)^2,  & a=0, \\
 e^{a(t+s_2)}(t+s_1+\tau) &\\
+\frac{b}{a^2}\big(e^{a(t+s_2)}-1-a(t+s_2)\big), & a\neq 0;
         \end{array}
       \right.
\]
and if $-s_1< t\le \tau-s_2$, then for  $a=0$ we obtain
\begin{align*}
R_t(s_1,s_2)&=\tau +\frac{b}{2}(t+s_2)^2 + \frac{b}{2}(t+s_1)^2\\
&\quad +\frac{b^2}{3}(t+s_1)^3+\frac{b^2}{2}(t+s_1)^2(s_2-s_1)
\end{align*}
while for   $a\neq 0$ we have
\begin{align*}
 R_t(s_1,s_2)&= e^{a(2t+s_1+s_2)}\tau +\frac{b}{a^2}e^{a(t+s_1)}\big(e^{a(t+s_2)}-1-a(t+s_2)\big)\\
&\quad +\frac{b}{a^2}\big(e^{a(t+s_2)}-\frac{b}{a}\big)\big(e^{a(t+s_1)}-1-a(t+s_1)\big)\\
&\quad + \frac{b^2}{2a^3}e^{a(s_2-s_1)}\big(e^{a(t+s_1)}-1\big)^2.
\end{align*}
In particular, we see that for $t\in [0,\tau]$
\[
R_t(-\tau,0)=\left\{
         \begin{array}{ll}
          t +\frac{b}{2}t^2,  & a=0, \\
 e^{at}t +\frac{b}{a^2}\big(e^{at}-1-at\big), & a\neq 0,
         \end{array}
       \right.
\]
and
\[
\sigma^2(t)=\left\{
         \begin{array}{ll}
          \tau +bt^2+\frac{b^2}{3}t^3,  & a=0, \\
 e^{2at}\tau +\frac{2b}{a^2}e^{at}\big(e^{at}-1-at\big)&\\
+\frac{b^2}{2a^3}\big((e^{at}-1)^2-2(e^{at}-1-at)\big), & a\neq 0.
         \end{array}
       \right.
\]
Since $\sigma^2(t)>0$ and $\det(Q_t)>0$ for $t\in (0,\tau]$, the densities $f$ in \eqref{e:densf} and $f_\nu$ in \eqref{e:densfn} are well defined. However, it should be noted that the Frobenius-Perron operator as in~\eqref{eq.fpdde1} will not be well defined if we take as a reference measure $\me$ on $C$ the Gaussian measure $\mu_0$. Suppose, on the contrary,  that $P^t$ is well defined.
According to the Feldman-H\'ajek theorem\cite{bogachev98},  two Gaussian measures are either equivalent (mutually absolutely continuous) or singular.
The Gaussian measure $\mu_t$, being the distribution of the process $S_t\xi$ with covariance function $R_t$, is equivalent to the Gaussian measure $\mu_0$\cite{shepp66,park72} if and only if
\[
R_t(s_1,s_2)=\min\{s_1,s_2\}+\tau -\int_{-\tau}^{s_1}\int_{-\tau}^{s_2}  K_t(r_1,r_2)dr_1dr_2,
\]
where $K_t\in L^2([-\tau,0]^2)$ is a symmetric function such that the integral operator
\[
K_t\phi(s)=\int_{-\tau}^0 K_t(s,r)\phi(r)dr
\]
does not have eigenvalue $1$. The kernel $K_t$ is unique and  for a.e. $(s_1,s_2)$ satisfies
\[
K_t(s_1,s_2)=-\frac{\partial}{\partial s_1}\frac{\partial}{\partial s_2}R_t(s_1,s_2).
\]
The formula that we have obtained for $R_t$ implies that the function $K_t$ does not exist in this example, thus showing that $\mu_t$ is singular with respect to $\me=\mu_0$. Hence, there exists a measurable set $A\subset C$ such that $\me(A)=0$ and $\me(S_t^{-1}(A))=1$,  contradicting~\eqref{eq.fpdde1} with $f(\phi)\equiv 1$. Consequently, the approach presented here in Section \ref{s:weak} may be more effective than the ones we briefly mentioned in Section~\ref{sec:formulate}.

\section{Summary}\label{sec:summary}

Here we have highlighted an open mathematical problem.  Namely, how can one formulate and study the evolution of densities in systems with dynamics that contain delays.  This is not simply an abstract mathematical problem devoid of interest.  Rather it is of prime scientific interest because of the increasing prevalence of studies of systems whose dynamics contain significant time delays, and the fact that we have no way of theoretically treating these systems.  Various approximations are available, and we have briefly discussed some of these while mentioning that there is a much more extensive consideration of these\cite{losson2020}.  We have also presented a tentative new approach to the problem in Section \ref{s:weak}.

Tied into this problem is the related issue of the lack of a well developed theory of bifurcation patterns of {\it densities} analogous to what exists for trajectories of dynamical systems.  In this vein we have also raised the question of whether it is possible for semi-dynamical systems to display a {\it chaotic} pattern of density evolution--currently an open question.

\section{Data Availability Statement.}  Data sharing is not applicable to this article as no new data were created or analyzed in this study.

\begin{acknowledgments}
We would like to acknowledge  research support from  the NSERC (Natural Sciences and Engineering Research Council of Canada), MITACS (Mathematics of Information Technology and Complex Systems),  the Alexander von Humboldt Stiftung and the National Science Centre (Poland) with grant No. 2017/27/B/ST1/00100.  Additionally we have greatly benefited from conversations with our colleagues Jinzhi Lei, J\'{e}r\^{o}me Losson, Nicholas Provatas, and Richard Taylor.  Finally the comments of two anonymous referees were most helpful and gratefully received.
\end{acknowledgments}

%


\begin{thebibliography}{42}%
\makeatletter
\providecommand \@ifxundefined [1]{%
 \@ifx{#1\undefined}
}%
\providecommand \@ifnum [1]{%
 \ifnum #1\expandafter \@firstoftwo
 \else \expandafter \@secondoftwo
 \fi
}%
\providecommand \@ifx [1]{%
 \ifx #1\expandafter \@firstoftwo
 \else \expandafter \@secondoftwo
 \fi
}%
\providecommand \natexlab [1]{#1}%
\providecommand \enquote  [1]{``#1''}%
\providecommand \bibnamefont  [1]{#1}%
\providecommand \bibfnamefont [1]{#1}%
\providecommand \citenamefont [1]{#1}%
\providecommand \href@noop [0]{\@secondoftwo}%
\providecommand \href [0]{\begingroup \@sanitize@url \@href}%
\providecommand \@href[1]{\@@startlink{#1}\@@href}%
\providecommand \@@href[1]{\endgroup#1\@@endlink}%
\providecommand \@sanitize@url [0]{\catcode `\\12\catcode `\$12\catcode
  `\&12\catcode `\#12\catcode `\^12\catcode `\_12\catcode `\%12\relax}%
\providecommand \@@startlink[1]{}%
\providecommand \@@endlink[0]{}%
\providecommand \url  [0]{\begingroup\@sanitize@url \@url }%
\providecommand \@url [1]{\endgroup\@href {#1}{\urlprefix }}%
\providecommand \urlprefix  [0]{URL }%
\providecommand \Eprint [0]{\href }%
\providecommand \doibase [0]{https://doi.org/}%
\providecommand \selectlanguage [0]{\@gobble}%
\providecommand \bibinfo  [0]{\@secondoftwo}%
\providecommand \bibfield  [0]{\@secondoftwo}%
\providecommand \translation [1]{[#1]}%
\providecommand \BibitemOpen [0]{}%
\providecommand \bibitemStop [0]{}%
\providecommand \bibitemNoStop [0]{.\EOS\space}%
\providecommand \EOS [0]{\spacefactor3000\relax}%
\providecommand \BibitemShut  [1]{\csname bibitem#1\endcsname}%
\let\auto@bib@innerbib\@empty
\bibitem [{\citenamefont {Hunt}\ and\ \citenamefont {Ott}(2015)}]{hunt2015}%
  \BibitemOpen
  \bibfield  {author} {\bibinfo {author} {\bibfnamefont {B.}~\bibnamefont
  {Hunt}}\ and\ \bibinfo {author} {\bibfnamefont {E.}~\bibnamefont {Ott}},\
  }\bibfield  {title} {\enquote {\bibinfo {title} {Defining chaos},}\
  }\href@noop {} {\bibfield  {journal} {\bibinfo  {journal} {Chaos}\ }\textbf
  {\bibinfo {volume} {5}},\ \bibinfo {pages} {097618} (\bibinfo {year}
  {2015})}\BibitemShut {NoStop}%
\bibitem [{\citenamefont {Lasota}\ and\ \citenamefont {Mackey}(1994)}]{LM94}%
  \BibitemOpen
  \bibfield  {author} {\bibinfo {author} {\bibfnamefont {A.}~\bibnamefont
  {Lasota}}\ and\ \bibinfo {author} {\bibfnamefont {M.~C.}\ \bibnamefont
  {Mackey}},\ }\href@noop {} {\emph {\bibinfo {title} {Chaos, fractals, and
  noise: Stochastic aspects of dynamics}}},\ \bibinfo {series} {Applied
  Mathematical Sciences}, Vol.~\bibinfo {volume} {97}\ (\bibinfo  {publisher}
  {Springer-Verlag},\ \bibinfo {address} {New York},\ \bibinfo {year}
  {1994})\BibitemShut {NoStop}%
\bibitem [{\citenamefont {Gardiner}(1983)}]{gardinerhandbook}%
  \BibitemOpen
  \bibfield  {author} {\bibinfo {author} {\bibfnamefont {C.}~\bibnamefont
  {Gardiner}},\ }\href@noop {} {\emph {\bibinfo {title} {Handbook of Stochastic
  Methods}}}\ (\bibinfo  {publisher} {Springer Verlag},\ \bibinfo {address}
  {Berlin, Heidelberg},\ \bibinfo {year} {1983})\BibitemShut {NoStop}%
\bibitem [{\citenamefont {Mallet-Paret}\ and\ \citenamefont
  {Nussbaum}(1986{\natexlab{a}})}]{mallet1986global}%
  \BibitemOpen
  \bibfield  {author} {\bibinfo {author} {\bibfnamefont {J.}~\bibnamefont
  {Mallet-Paret}}\ and\ \bibinfo {author} {\bibfnamefont {R.~D.}\ \bibnamefont
  {Nussbaum}},\ }\bibfield  {title} {\enquote {\bibinfo {title} {Global
  continuation and asymptotic behaviour for periodic solutions of a
  differential-delay equation},}\ }\href@noop {} {\bibfield  {journal}
  {\bibinfo  {journal} {Ann. Mat. Pura Appl. }\ }\textbf
  {\bibinfo {volume} {145}},\ \bibinfo {pages} {33--128} (\bibinfo {year}
  {1986}{\natexlab{a}})}\BibitemShut {NoStop}%
\bibitem [{\citenamefont {Sharkovski{\u\i}}, \citenamefont {Ma{\u\i}strenko},\
  and\ \citenamefont {Romanenko}(1986)}]{sharkovski-bk}%
  \BibitemOpen
  \bibfield  {author} {\bibinfo {author} {\bibfnamefont {A.~N.}\ \bibnamefont
  {Sharkovski{\u\i}}}, \bibinfo {author} {\bibfnamefont {Y.~L.}\ \bibnamefont
  {Ma{\u\i}strenko}},\ and\ \bibinfo {author} {\bibfnamefont {E.~Y.}\
  \bibnamefont {Romanenko}},\ }\href@noop {} {\emph {\bibinfo {title}
  {Raznostnye uravneniya i ikh prilozheniya}}}\ (\bibinfo  {publisher}
  {``Naukova Dumka''},\ \bibinfo {address} {Kiev},\ \bibinfo {year}
  {1986})\BibitemShut {NoStop}%
\bibitem [{\citenamefont {Mallet-Paret}\ and\ \citenamefont
  {Nussbaum}(1986{\natexlab{b}})}]{mallet1986bifurcation}%
  \BibitemOpen
  \bibfield  {author} {\bibinfo {author} {\bibfnamefont {J.}~\bibnamefont
  {Mallet-Paret}}\ and\ \bibinfo {author} {\bibfnamefont {R.~D.}\ \bibnamefont
  {Nussbaum}},\ }\bibfield  {title} {\enquote {\bibinfo {title} {A bifurcation
  gap for a singularly perturbed delay equation},}\ }in\ \href@noop {} {\emph
  {\bibinfo {booktitle} {Chaotic Dynamics and Fractals}}}\ (\bibinfo
  {publisher} {Elsevier},\ \bibinfo {year} {1986})\ pp.\ \bibinfo {pages}
  {263--286}\BibitemShut {NoStop}%
\bibitem [{\citenamefont {Kuznetsov}(2013)}]{kuznetsov2013elements}%
  \BibitemOpen
  \bibfield  {author} {\bibinfo {author} {\bibfnamefont {Y.~A.}\ \bibnamefont
  {Kuznetsov}},\ }\href@noop {} {\emph {\bibinfo {title} {Elements of applied
  bifurcation theory}}},\ Vol.\ \bibinfo {volume} {112}\ (\bibinfo  {publisher}
  {Springer Science \& Business Media},\ \bibinfo {year} {2013})\BibitemShut
  {NoStop}%
\bibitem [{Note1()}]{Note1}%
  \BibitemOpen
  \bibinfo {note} {$A$ is invariant if $S^{-1}(A)=A$}\BibitemShut {NoStop}%
\bibitem [{\citenamefont {Ito}, \citenamefont {Tanaka},\ and\ \citenamefont
  {Nakada}(1979{\natexlab{a}})}]{ito79a}%
  \BibitemOpen
  \bibfield  {author} {\bibinfo {author} {\bibfnamefont {S.}~\bibnamefont
  {Ito}}, \bibinfo {author} {\bibfnamefont {S.}~\bibnamefont {Tanaka}},\ and\
  \bibinfo {author} {\bibfnamefont {H.}~\bibnamefont {Nakada}},\ }\bibfield
  {title} {\enquote {\bibinfo {title} {On unimodal linear transformations and
  chaos. {I}},}\ }\href@noop {} {\bibfield  {journal} {\bibinfo  {journal}
  {Tokyo J. Math.}\ }\textbf {\bibinfo {volume} {2}},\ \bibinfo {pages}
  {221--239} (\bibinfo {year} {1979}{\natexlab{a}})}\BibitemShut {NoStop}%
\bibitem [{\citenamefont {Ito}, \citenamefont {Tanaka},\ and\ \citenamefont
  {Nakada}(1979{\natexlab{b}})}]{ito79b}%
  \BibitemOpen
  \bibfield  {author} {\bibinfo {author} {\bibfnamefont {S.}~\bibnamefont
  {Ito}}, \bibinfo {author} {\bibfnamefont {S.}~\bibnamefont {Tanaka}},\ and\
  \bibinfo {author} {\bibfnamefont {H.}~\bibnamefont {Nakada}},\ }\bibfield
  {title} {\enquote {\bibinfo {title} {On unimodal linear transformations and
  chaos. {II}},}\ }\href@noop {} {\bibfield  {journal} {\bibinfo  {journal}
  {Tokyo J. Math.}\ }\textbf {\bibinfo {volume} {2}},\ \bibinfo {pages}
  {241--259} (\bibinfo {year} {1979}{\natexlab{b}})}\BibitemShut {NoStop}%
\bibitem [{\citenamefont {Yoshida}, \citenamefont {Mori},\ and\ \citenamefont
  {Shigematsu}(1983)}]{yoshida83}%
  \BibitemOpen
  \bibfield  {author} {\bibinfo {author} {\bibfnamefont {T.}~\bibnamefont
  {Yoshida}}, \bibinfo {author} {\bibfnamefont {H.}~\bibnamefont {Mori}},\ and\
  \bibinfo {author} {\bibfnamefont {H.}~\bibnamefont {Shigematsu}},\ }\bibfield
   {title} {\enquote {\bibinfo {title} {Analytic study of chaos of the tent
  map: band structures, power spectra, and critical behaviors},}\ }\href@noop
  {} {\bibfield  {journal} {\bibinfo  {journal} {J. Statist. Phys.}\ }\textbf
  {\bibinfo {volume} {31}},\ \bibinfo {pages} {279--308} (\bibinfo {year}
  {1983})}\BibitemShut {NoStop}%
\bibitem [{\citenamefont {Provatas}\ and\ \citenamefont
  {Mackey}(1991{\natexlab{a}})}]{provatas91a}%
  \BibitemOpen
  \bibfield  {author} {\bibinfo {author} {\bibfnamefont {N.}~\bibnamefont
  {Provatas}}\ and\ \bibinfo {author} {\bibfnamefont {M.~C.}\ \bibnamefont
  {Mackey}},\ }\bibfield  {title} {\enquote {\bibinfo {title} {Asymptotic
  periodicity and banded chaos},}\ }\href@noop {} {\bibfield  {journal}
  {\bibinfo  {journal} {Phys. D}\ }\textbf {\bibinfo {volume} {53}},\ \bibinfo
  {pages} {295--318} (\bibinfo {year} {1991}{\natexlab{a}})}\BibitemShut
  {NoStop}%
\bibitem [{\citenamefont {Arnold}(1998)}]{arnold98}%
  \BibitemOpen
  \bibfield  {author} {\bibinfo {author} {\bibfnamefont {L.}~\bibnamefont
  {Arnold}},\ }\href@noop {} {\emph {\bibinfo {title} {Random dynamical
  systems}}},\ Springer Monographs in Mathematics\ (\bibinfo  {publisher}
  {Springer-Verlag},\ \bibinfo {address} {Berlin},\ \bibinfo {year}
  {1998})\BibitemShut {NoStop}%
\bibitem [{\citenamefont {Mackey}(2009)}]{mackey2009exploring}%
  \BibitemOpen
  \bibfield  {author} {\bibinfo {author} {\bibfnamefont {M.~C.}\ \bibnamefont
  {Mackey}},\ }\bibfield  {title} {\enquote {\bibinfo {title} {Exploring the
  world with mathematics},}\ }\href@noop {} {\bibfield  {journal} {\bibinfo
  {journal} {Ann. Math. Sil.}\ }\textbf {\bibinfo {volume} {23}},\ \bibinfo
  {pages} {11--42} (\bibinfo {year} {2009})}\BibitemShut {NoStop}%
\bibitem [{\citenamefont {Mackey}\ \emph {et~al.}(2012)\citenamefont {Mackey},
  \citenamefont {Tyran-Kami{\'n}ska}, \citenamefont {Walther} \emph
  {et~al.}}]{mackey2012mathematical}%
  \BibitemOpen
  \bibfield  {author} {\bibinfo {author} {\bibfnamefont {M.~C.}\ \bibnamefont
  {Mackey}}, \bibinfo {author} {\bibfnamefont {M.}~\bibnamefont
  {Tyran-Kami{\'n}ska}}, \bibinfo {author} {\bibfnamefont {H.-O.}\ \bibnamefont
  {Walther}}, \emph {et~al.},\ }\bibfield  {title} {\enquote {\bibinfo {title}
  {The mathematical legacy of {A}ndrzej {L}asota},}\ }\href@noop {} {\bibfield
  {journal} {\bibinfo  {journal} {Wiad. Mat.}\ }\textbf {\bibinfo {volume}
  {48}},\ \bibinfo {pages} {143} (\bibinfo {year} {2012})}\BibitemShut
  {NoStop}%
\bibitem [{\citenamefont {Mackey}(2016)}]{mackey2016adventures}%
  \BibitemOpen
  \bibfield  {author} {\bibinfo {author} {\bibfnamefont {M.~C.}\ \bibnamefont
  {Mackey}},\ }\bibfield  {title} {\enquote {\bibinfo {title} {Adventures in
  {P}oland: Having fun and doing research with {A}ndrzej {L}asota},}\
  }\href@noop {} {\bibfield  {journal} {\bibinfo  {journal} {Mat. Appl.
  (Warsaw)}\ }\textbf {\bibinfo {volume} {35}},\ \bibinfo {pages} {5--32}
  (\bibinfo {year} {2016})}\BibitemShut {NoStop}%
\bibitem [{\citenamefont {Losson}\ and\ \citenamefont
  {Mackey}(1995{\natexlab{a}})}]{LM95}%
  \BibitemOpen
  \bibfield  {author} {\bibinfo {author} {\bibfnamefont {J.}~\bibnamefont
  {Losson}}\ and\ \bibinfo {author} {\bibfnamefont {M.~C.}\ \bibnamefont
  {Mackey}},\ }\bibfield  {title} {\enquote {\bibinfo {title} {Coupled map
  lattices as models of deterministic and stochastic differential delay
  equations},}\ }\href@noop {} {\bibfield  {journal} {\bibinfo  {journal}
  {Phys. Rev. E}\ }\textbf {\bibinfo {volume} {52}},\ \bibinfo {pages}
  {115--128} (\bibinfo {year} {1995}{\natexlab{a}})}\BibitemShut {NoStop}%
\bibitem [{\citenamefont {Lasota}\ and\ \citenamefont
  {Mackey}(1987)}]{almcm87}%
  \BibitemOpen
  \bibfield  {author} {\bibinfo {author} {\bibfnamefont {A.}~\bibnamefont
  {Lasota}}\ and\ \bibinfo {author} {\bibfnamefont {M.~C.}\ \bibnamefont
  {Mackey}},\ }\bibfield  {title} {\enquote {\bibinfo {title} {Noise and
  statistical periodicity},}\ }\href@noop {} {\bibfield  {journal} {\bibinfo
  {journal} {Phys. D}\ }\textbf {\bibinfo {volume} {28}},\ \bibinfo {pages}
  {143--154} (\bibinfo {year} {1987})}\BibitemShut {NoStop}%
\bibitem [{\citenamefont {Provatas}\ and\ \citenamefont
  {Mackey}(1991{\natexlab{b}})}]{provatas91b}%
  \BibitemOpen
  \bibfield  {author} {\bibinfo {author} {\bibfnamefont {N.}~\bibnamefont
  {Provatas}}\ and\ \bibinfo {author} {\bibfnamefont {M.~C.}\ \bibnamefont
  {Mackey}},\ }\bibfield  {title} {\enquote {\bibinfo {title} {Noise-induced
  asymptotic periodicity in a piecewise linear map},}\ }\href@noop {}
  {\bibfield  {journal} {\bibinfo  {journal} {J. Statist. Phys.}\ }\textbf
  {\bibinfo {volume} {63}},\ \bibinfo {pages} {585--612} (\bibinfo {year}
  {1991}{\natexlab{b}})}\BibitemShut {NoStop}%
\bibitem [{\citenamefont {Lei}\ and\ \citenamefont
  {Mackey}(2011)}]{lei2011deterministic}%
  \BibitemOpen
  \bibfield  {author} {\bibinfo {author} {\bibfnamefont {J.}~\bibnamefont
  {Lei}}\ and\ \bibinfo {author} {\bibfnamefont {M.~C.}\ \bibnamefont
  {Mackey}},\ }\bibfield  {title} {\enquote {\bibinfo {title} {Deterministic
  {B}rownian motion generated from differential delay equations},}\ }\href@noop
  {} {\bibfield  {journal} {\bibinfo  {journal} {Phys. Rev. E}\ }\textbf
  {\bibinfo {volume} {84}},\ \bibinfo {pages} {041105} (\bibinfo {year}
  {2011})}\BibitemShut {NoStop}%
\bibitem [{\citenamefont {Mackey}\ and\ \citenamefont
  {Tyran-Kami{\'n}ska}(2006)}]{mackeytyran}%
  \BibitemOpen
  \bibfield  {author} {\bibinfo {author} {\bibfnamefont {M.~C.}\ \bibnamefont
  {Mackey}}\ and\ \bibinfo {author} {\bibfnamefont {M.}~\bibnamefont
  {Tyran-Kami{\'n}ska}},\ }\bibfield  {title} {\enquote {\bibinfo {title}
  {Deterministic {B}rownian motion: {T}he effects of perturbing a dynamical
  system by a chaotic semi-dynamical system},}\ }\href@noop {} {\bibfield
  {journal} {\bibinfo  {journal} {Phys. Rep.}\ }\textbf {\bibinfo {volume}
  {422}},\ \bibinfo {pages} {167--222} (\bibinfo {year} {2006})}\BibitemShut
  {NoStop}%
\bibitem [{\citenamefont {Tyran-Kami\'{n}ska}(2014)}]{mtyran14}%
  \BibitemOpen
  \bibfield  {author} {\bibinfo {author} {\bibfnamefont {M.}~\bibnamefont
  {Tyran-Kami\'{n}ska}},\ }\bibfield  {title} {\enquote {\bibinfo {title}
  {Diffusion and deterministic systems},}\ }\href
  {https://doi.org/10.1051/mmnp/20149110} {\bibfield  {journal} {\bibinfo
  {journal} {Math. Model. Nat. Phenom.}\ }\textbf {\bibinfo {volume} {9}},\
  \bibinfo {pages} {139--150} (\bibinfo {year} {2014})}\BibitemShut {NoStop}%
\bibitem [{\citenamefont {Losson}\ \emph {et~al.}(2020)\citenamefont {Losson},
  \citenamefont {Mackey}, \citenamefont {Taylor},\ and\ \citenamefont
  {Tyran-Kami{\'n}ska}}]{losson2020}%
  \BibitemOpen
  \bibfield  {author} {\bibinfo {author} {\bibfnamefont {J.}~\bibnamefont
  {Losson}}, \bibinfo {author} {\bibfnamefont {M.~C.}\ \bibnamefont {Mackey}},
  \bibinfo {author} {\bibfnamefont {R.}~\bibnamefont {Taylor}},\ and\ \bibinfo
  {author} {\bibfnamefont {M.}~\bibnamefont {Tyran-Kami{\'n}ska}},\ }\href@noop
  {} {\emph {\bibinfo {title} {Density evolution under delayed dynamics: {A}n
  open problem}}}\ (\bibinfo  {publisher} {Springer Verlag},\ \bibinfo {year}
  {2020})\BibitemShut {NoStop}%
\bibitem [{\citenamefont {Driver}(1977)}]{driver2012ordinary}%
  \BibitemOpen
  \bibfield  {author} {\bibinfo {author} {\bibfnamefont {R.~D.}\ \bibnamefont
  {Driver}},\ }\href@noop {} {\emph {\bibinfo {title} {Ordinary and delay
  differential equations}}},\ Vol.~\bibinfo {volume} {20}\ (\bibinfo
  {publisher} {Springer Science \& Business Media},\ \bibinfo {year}
  {1977})\BibitemShut {NoStop}%
\bibitem [{\citenamefont {Hopf}(1952)}]{hopf1952statistical}%
  \BibitemOpen
  \bibfield  {author} {\bibinfo {author} {\bibfnamefont {E.}~\bibnamefont
  {Hopf}},\ }\bibfield  {title} {\enquote {\bibinfo {title} {Statistical
  hydromechanics and functional calculus},}\ }\href@noop {} {\bibfield
  {journal} {\bibinfo  {journal} {J. Rat. Mech. Anal.}\ }\textbf {\bibinfo
  {volume} {1}},\ \bibinfo {pages} {87--123} (\bibinfo {year}
  {1952})}\BibitemShut {NoStop}%
\bibitem [{\citenamefont {Losson}\ and\ \citenamefont
  {Mackey}(1992)}]{losson92}%
  \BibitemOpen
  \bibfield  {author} {\bibinfo {author} {\bibfnamefont {J.}~\bibnamefont
  {Losson}}\ and\ \bibinfo {author} {\bibfnamefont {M.~C.}\ \bibnamefont
  {Mackey}},\ }\bibfield  {title} {\enquote {\bibinfo {title} {A {H}opf-like
  equation and perturbation theory for differential delay equations},}\
  }\href@noop {} {\bibfield  {journal} {\bibinfo  {journal} {J. Statist.
  Phys.}\ }\textbf {\bibinfo {volume} {69}},\ \bibinfo {pages} {1025--1046}
  (\bibinfo {year} {1992})}\BibitemShut {NoStop}%
\bibitem [{\citenamefont {Kac}(1980)}]{Kac}%
  \BibitemOpen
  \bibfield  {author} {\bibinfo {author} {\bibfnamefont {M.}~\bibnamefont
  {Kac}},\ }\href@noop {} {\emph {\bibinfo {title} {Integration in function
  spaces and some of its applications}}}\ (\bibinfo  {publisher} {Accademia
  Nazionale dei Lincei, Scuola Normale Superiore},\ \bibinfo {address} {Pisa},\
  \bibinfo {year} {1980})\BibitemShut {NoStop}%
\bibitem [{\citenamefont {Brunovsky}\ and\ \citenamefont
  {Komornik}(1984)}]{BK84}%
  \BibitemOpen
  \bibfield  {author} {\bibinfo {author} {\bibfnamefont {P.}~\bibnamefont
  {Brunovsky}}\ and\ \bibinfo {author} {\bibfnamefont {J.}~\bibnamefont
  {Komornik}},\ }\bibfield  {title} {\enquote {\bibinfo {title} {Ergodicity and
  exactness of the shift on ${C}[0,\infty)$ and the semiflow of a first-order
  partial differential equation},}\ }\href@noop {} {\bibfield  {journal}
  {\bibinfo  {journal} {J.~Math. Anal. Appl.}\
  }\textbf {\bibinfo {volume} {104}},\ \bibinfo {pages} {235--245} (\bibinfo
  {year} {1984})}\BibitemShut {NoStop}%
\bibitem [{\citenamefont {Rudnicki}(1985)}]{Rud85}%
  \BibitemOpen
  \bibfield  {author} {\bibinfo {author} {\bibfnamefont {R.}~\bibnamefont
  {Rudnicki}},\ }\bibfield  {title} {\enquote {\bibinfo {title} {Invariant
  measures for the flow of a first order partial differential equation},}\
  }\href@noop {} {\bibfield  {journal} {\bibinfo  {journal} {Ergodic Theory Dynam. Systems}\ }\textbf {\bibinfo {volume} {5}},\ \bibinfo {pages}
  {437--443} (\bibinfo {year} {1985})}\BibitemShut {NoStop}%
\bibitem [{\citenamefont {Rudnicki}(1987)}]{Rud87}%
  \BibitemOpen
  \bibfield  {author} {\bibinfo {author} {\bibfnamefont {R.}~\bibnamefont
  {Rudnicki}},\ }\bibfield  {title} {\enquote {\bibinfo {title} {An abstract
  {W}iener measure invariant under a partial differential equation},}\
  }\href@noop {} {\bibfield  {journal} {\bibinfo  {journal} {Bull. Polish Acad. Sci. Math.}\ }\textbf {\bibinfo {volume} {35}},\
  \bibinfo {pages} {289--295} (\bibinfo {year} {1987})}\BibitemShut {NoStop}%
\bibitem [{\citenamefont {Rudnicki}(1988)}]{Rud88}%
  \BibitemOpen
  \bibfield  {author} {\bibinfo {author} {\bibfnamefont {R.}~\bibnamefont
  {Rudnicki}},\ }\bibfield  {title} {\enquote {\bibinfo {title} {Strong ergodic
  properties of a first-order partial differential equation},}\ }\href@noop {}
  {\bibfield  {journal} {\bibinfo  {journal} {J. Math. Anal. Appl.}\ }\textbf {\bibinfo {volume} {133}},\ \bibinfo {pages}
  {14--26} (\bibinfo {year} {1988})}\BibitemShut {NoStop}%
\bibitem [{\citenamefont {Losson}\ and\ \citenamefont
  {Mackey}(1995{\natexlab{b}})}]{losson1995evolution}%
  \BibitemOpen
  \bibfield  {author} {\bibinfo {author} {\bibfnamefont {J.}~\bibnamefont
  {Losson}}\ and\ \bibinfo {author} {\bibfnamefont {M.~C.}\ \bibnamefont
  {Mackey}},\ }\bibfield  {title} {\enquote {\bibinfo {title} {Evolution of
  probability densities in stochastic coupled map lattices},}\ }\href@noop {}
  {\bibfield  {journal} {\bibinfo  {journal} {Phys. Rev. E}\ }\textbf {\bibinfo
  {volume} {52}},\ \bibinfo {pages} {1403} (\bibinfo {year}
  {1995}{\natexlab{b}})}\BibitemShut {NoStop}%
\bibitem [{\citenamefont {Ionescu~Tulcea}\ and\ \citenamefont
  {Marinescu}(1950)}]{ionescu50}%
  \BibitemOpen
  \bibfield  {author} {\bibinfo {author} {\bibfnamefont {C.~T.}\ \bibnamefont
  {Ionescu~Tulcea}}\ and\ \bibinfo {author} {\bibfnamefont {G.}~\bibnamefont
  {Marinescu}},\ }\bibfield  {title} {\enquote {\bibinfo {title} {Th\'eorie
  ergodique pour des classes d'op\'erations non compl\`etement continues},}\
  }\href@noop {} {\bibfield  {journal} {\bibinfo  {journal} {Ann. of Math.
  (2)}\ }\textbf {\bibinfo {volume} {52}},\ \bibinfo {pages} {140--147}
  (\bibinfo {year} {1950})}\BibitemShut {NoStop}%
\bibitem [{\citenamefont {Dellnitz}\ and\ \citenamefont
  {Junge}(1999)}]{dellnitz1999approximation}%
  \BibitemOpen
  \bibfield  {author} {\bibinfo {author} {\bibfnamefont {M.}~\bibnamefont
  {Dellnitz}}\ and\ \bibinfo {author} {\bibfnamefont {O.}~\bibnamefont
  {Junge}},\ }\bibfield  {title} {\enquote {\bibinfo {title} {On the
  approximation of complicated dynamical behavior},}\ }\href@noop {} {\bibfield
   {journal} {\bibinfo  {journal} {SIAM J. Numer. Anal.}\
  }\textbf {\bibinfo {volume} {36}},\ \bibinfo {pages} {491--515} (\bibinfo
  {year} {1999})}\BibitemShut {NoStop}%
\bibitem [{\citenamefont {Dellnitz}, \citenamefont {Hessel-Von~Molo},\ and\
  \citenamefont {Ziessler}(2016)}]{dellnitz2015computation}%
  \BibitemOpen
  \bibfield  {author} {\bibinfo {author} {\bibfnamefont {M.}~\bibnamefont
  {Dellnitz}}, \bibinfo {author} {\bibfnamefont {M.}~\bibnamefont
  {Hessel-Von~Molo}},\ and\ \bibinfo {author} {\bibfnamefont {A.}~\bibnamefont
  {Ziessler}},\ }\bibfield  {title} {\enquote {\bibinfo {title} {On the
  computation of attractors for delay differential equations},}\ }\href@noop {}
  {\bibfield  {journal} {\bibinfo  {journal} {J. Comput. Dyn.}\ }\textbf
  {\bibinfo {volume} {3}},\ \bibinfo {pages} {93--112} (\bibinfo {year}
  {2016})}\BibitemShut {NoStop}%
\bibitem [{\citenamefont {Hale}\ and\ \citenamefont
  {Verduyn~Lunel}(1993)}]{hale-lunel}%
  \BibitemOpen
  \bibfield  {author} {\bibinfo {author} {\bibfnamefont {J.~K.}\ \bibnamefont
  {Hale}}\ and\ \bibinfo {author} {\bibfnamefont {S.~M.}\ \bibnamefont
  {Verduyn~Lunel}},\ }\href@noop {} {\emph {\bibinfo {title} {Introduction to
  functional-differential equations}}},\ \bibinfo {series} {Applied
  Mathematical Sciences}, Vol.~\bibinfo {volume} {99}\ (\bibinfo  {publisher}
  {Springer-Verlag},\ \bibinfo {address} {New York},\ \bibinfo {year}
  {1993})\BibitemShut {NoStop}%
\bibitem [{\citenamefont {Dynkin}(1965)}]{dynkin}%
  \BibitemOpen
  \bibfield  {author} {\bibinfo {author} {\bibfnamefont {E.~B.}\ \bibnamefont
  {Dynkin}},\ }\href@noop {} {\emph {\bibinfo {title} {Markov processes.
  {V}ols. {I}, {II}}}},\ \bibinfo {series} {Die Grundlehren der Mathematischen
  Wissenschaften, B\"{a}nde 121}, Vol.\ \bibinfo {volume} {122}\ (\bibinfo
  {publisher} {Academic Press Inc., Publishers, New York; Springer-Verlag,
  Berlin-G\"{o}ttingen-Heidelberg},\ \bibinfo {year} {1965})\BibitemShut
  {NoStop}%
\bibitem [{\citenamefont {Mohammed}(1984)}]{mohammed}%
  \BibitemOpen
  \bibfield  {author} {\bibinfo {author} {\bibfnamefont {S.~E.~A.}\
  \bibnamefont {Mohammed}},\ }\href@noop {} {\emph {\bibinfo {title}
  {Stochastic functional differential equations}}},\ \bibinfo {series}
  {Research Notes in Mathematics}, Vol.~\bibinfo {volume} {99}\ (\bibinfo
  {publisher} {Pitman Advanced Publishing Program, Boston, MA},\ \bibinfo
  {year} {1984})\BibitemShut {NoStop}%
\bibitem [{\citenamefont {Bogachev}(1998)}]{bogachev98}%
  \BibitemOpen
  \bibfield  {author} {\bibinfo {author} {\bibfnamefont {V.~I.}\ \bibnamefont
  {Bogachev}},\ }\href {https://doi.org/10.1090/surv/062} {\emph {\bibinfo
  {title} {Gaussian measures}}},\ \bibinfo {series} {Mathematical Surveys and
  Monographs}, Vol.~\bibinfo {volume} {62}\ (\bibinfo  {publisher} {American
  Mathematical Society, Providence, RI},\ \bibinfo {year} {1998})\BibitemShut
  {NoStop}%
\bibitem [{\citenamefont {Hayes}(1950)}]{hayes1950roots}%
  \BibitemOpen
  \bibfield  {author} {\bibinfo {author} {\bibfnamefont {N.}~\bibnamefont
  {Hayes}},\ }\bibfield  {title} {\enquote {\bibinfo {title} {Roots of the
  transcendental equation associated with a certain difference-differential
  equation},}\ }\href@noop {} {\bibfield  {journal} {\bibinfo  {journal} {J.
  London Math. Soc. (2)}\ }\textbf {\bibinfo {volume} {1}},\ \bibinfo {pages}
  {226--232} (\bibinfo {year} {1950})}\BibitemShut {NoStop}%
\bibitem [{\citenamefont {Shepp}(1966)}]{shepp66}%
  \BibitemOpen
  \bibfield  {author} {\bibinfo {author} {\bibfnamefont {L.~A.}\ \bibnamefont
  {Shepp}},\ }\bibfield  {title} {\enquote {\bibinfo {title}
  {Radon-{N}ikod\'{y}m derivatives of {G}aussian measures},}\ }\href
  {https://doi.org/10.1214/aoms/1177699516} {\bibfield  {journal} {\bibinfo
  {journal} {Ann. Math. Statist.}\ }\textbf {\bibinfo {volume} {37}},\ \bibinfo
  {pages} {321--354} (\bibinfo {year} {1966})}\BibitemShut {NoStop}%
\bibitem [{\citenamefont {Park}(1972)}]{park72}%
  \BibitemOpen
  \bibfield  {author} {\bibinfo {author} {\bibfnamefont {W.~J.}\ \bibnamefont
  {Park}},\ }\bibfield  {title} {\enquote {\bibinfo {title} {On the equivalence
  of {G}aussian processes with factorable covariance functions},}\ }\href
  {https://doi.org/10.2307/2038345} {\bibfield  {journal} {\bibinfo  {journal}
  {Proc. Amer. Math. Soc.}\ }\textbf {\bibinfo {volume} {32}},\ \bibinfo
  {pages} {275--279} (\bibinfo {year} {1972})}\BibitemShut {NoStop}%
\end{thebibliography}

\end{document}